%% file: tropDaulity_arx.tex
\documentclass[11pt]{article}
\input tropMacro.tex

\title{  \begin{center} Duality of Tropical Curves \end{center}}

\author{Zur Izhakian\thanks{The School of Mathematical \& Computer Science, Raymond
and Beverly Sackler Faculty of Exact Sciences, Tel Aviv
University, Ramat Aviv, Tel Aviv, 69978, Israel. Email:
zzur@post.tau.ac.il.}
\thanks{The author has been supported by The German-Israeli Foundation for Research and
Development by Hermann Minkowski Minerva Center for Geometry at
Tel Aviv University.}}
\date{}
\begin{document}
\maketitle

%******************************* abstract *********************************
\begin{abstract}
Duality of curves is an important aspect of the ``classical''
algebraic geometry. In this paper, using this foundation, the
duality of tropical polynomials is constructed to introduce the
duality of Non-Archimedean curves. Using the development of an
algebraic ``mechanism'', based on ``distortion'' values, geometric
and convexity properties are analyzed. Specifically, we discuss
some significant aspects referring to quadrics with respect to
their dual objects. This topic also includes the induced dual
subdivision of the corresponding Newton Polytope and its
compatible properties. Finally, a regularity of tropical curves in
the duality sense is generally defined and studied for families of
tropical quadrics.
\end{abstract}

%******************************* keywords *********************************
%\emph{Remarks}: \\....
%******************************* keywords *********************************
%\begin{keywords}
%\emph{{\bf Keywords}}: \\....
%\end{keywords}

%******************************* AMS *********************************
%\begin{AMS}
%\emph{{\bf AMS}} : \\....
%\end{AMS}

%******************************* section ******************************
%******************************* BODY *********************************
%******************************* SECT 1 *********************************
\secSpc
\section*{Instruction}
Over the last years, intensive development in the studies of the
tropical algebraic geometry has been made. The special ``nature''
of the objects belonging to this geometry enables progress in
varied directions and by using different concepts. The main
significance of these tropical entities is their being geometric
``images'' of algebraic objects and they concurrently comprise
combinatorial attributes. Recently, much efforts have been
invested to characterize the tropical analogous to ``classical''
results and to determine the various connections between these two
``worlds''. The purpose of this paper is to introduce the tropical
analogous to the ``classical'' curves' duality in algebraic
geometry, and to outline the connection between these dualities.

\parSpc
As known, duality between algebraic objects is a powerful tool
with many applications. This is the motivation of our work --
adapting this method to the tropical case as well. Thus,
understanding the geometric linkage between algebraic and tropical
objects, especially between dual objects, is a preliminary step
toward developing a geometric duality over polyhedra. This duality
should also agree with the algebraic duality.

\parSpc
Algebraic objects are formally, elements of the geometry over the
fields ($\Fld,+,\cdot$), where the  tropical ones are those over
the geometry of the semi-ring ($\Real,max,+$) -- the semi-ring
which contains the {\em Max-Plus Algebra}
\cite{Cuninghame2004,Sturmfels6366}. The fundamental objects in
this geometry are {\em Polyhedral Complexes}, where their behavior
resembles the complex algebraic varieties
\cite{Mikhalkin8225,Shustin1278}. Moreover, they may be concerned
as the ``images'' of the \emph{Non-Archimedean} valuation of some
``superior'' algebraic varieties. Historically, these objects
carry the name \emph{Non-Archimedean Amebas}. Using this approach,
a direct relationship exits between the following involved
components, the tropical hyper-surface $V(\tF)$ the corresponding
\emph{Newton Polytope} $\Polygon_\tF$ and the induced subdivision
$\tS_\tF$ over $\Polygon_\tF$ \cite{Mikhalkin01,Passare00}.
Namely, in this theory, algebraic and combinatorial considerations
together with geometric observations are all composed together.

\parSpc
This area of research was formally introduced for complex amebas
in 1994 by Gelfand, Kapranov and Zelevinsky \cite{Gelfand94}.
Later Kapranov firstly presented the notion of Non-Archimedean
Amebas \cite{Kapranov2000} which may be understood as the
``spines'' of the complex amebas \cite{Mikhalkin8225}. Since then,
miscellaneous aspects within this topic have been dealt within
several papers. Sturmfels, Speyer and Develin studied matters of
tropical algebra \cite{Develin8254,Speyer4218,Sturmfels02}.
Enumerative geometry has been dealt by Mikhalkin
\cite{Mikhalkin05,Mikhalkin03}. Shustin discussed the patchworking
of Non-Arachnidan amoebas \cite{Shustin1278} in the algebraic
manner, while Itenberg regarded this issues combinatorially
\cite{Itenberg96}. These are few of the many directions which
recently concerned researchers and are spread along varied fields
of study.

\parSpc
The goal of this paper is to introduce the geometric duality of
tropical curves and the conceptional connection between this
duality and that which is applied to algebraic objects. Since
tropical objects are basically piece-wise linear, aiming to
develop a duality we join together algebraical and combinatorial
methods. Despite the fact that tropical notion of tangency has not
been phrased properly yet, by relying on the duality of complex
algebraic varieties and their linkage to tropical objects we
overcome this obstacle and define the duality for tropical
polynomials. That is our point of departure for the discussion in
this paper.

\parSpc
Specifically, for a given tropical curve, $V(\tF)$ the properties
of the Non-Archimedean valuation are used to ``produce''
compatible algebraic varieties and so dose their dual, and these
are ``translated backward'' to attain the tropical dual objects
described by $V(\dtF)$. Eventually, we show that this whole
procedure is reduced to pure computation in terms of tropical
operations. This approach disregards convexity considerations and
yields only the pre-tropical duality (i.e. the dual ``tropical
polynomial'', $\dtF$). Thus, in order to obtain the proper
tropical varieties we should involve geometric considerations
which include convexity properties.

\parSpc
In order to observe the convexity properties of a given object we
develop an algebraic ``mechanism'', expressed in the sense of
``distortion'' values, which serves to refine the required
convexity attributes. Moreover, this ``mechanism'' appears to be
useful for additional geometric analysis. Using these methods, we
study the geometric duality of tropical curves via the
 subdivisions of their corresponding Newton polytopes
(i.e. the dual subdivision $S_{\dtF}$). This will be done mainly
by focusing on ``interesting'' families of subdivisions which are
significant for other researches. The speciality of these
subdivisions is that they can be described in the terms of the
polytope's nodes, and this makes the convexity analysis more
convenient.

\parSpc
Explicitly for the case of quadrics, we prove that a subdivision
which is induced by a tropical polynomial has a one to one
correspondence with the tropical curve. Moreover
\begin{itemize}
  \item The dual subdivision $S_{\dtF}$ of a  subdivision $S_{\tF}$  that is \full
  \; (\emty)
  is also \full \; (\emty).

%  \item The dual subdivision $S_{\dtF}$ of a  subdivision that is \emty \;
%  $S_{\tF}$ is also \emty.
\end{itemize}
These types of subdivisions preserve their preliminary properties
of convexity under the dual transformation. Advancing these ideas
and enforcing the additional requirement in which not only that
the resulting subdivision is regarded but also its specifier, we
characterize the regularity of tropical curves in the sense of
duality. This notion refers to general curves and it is discussed
in detail for the case of  quadrics.

\parSpc
Using the above notion, the regularity occurs when we apply the
duality twice for a curve, and the difference between the
corresponding specifiers of the result and those of the primal
curve is a constant. Specifically we show that
\begin{itemize}
  \item a curve that corresponds to a subdivision which is \emty \; is
  never regular while,
  \item by enforcing additional restrictions, a subdivision which is \full \;
   can correspond to a regular curve
\end{itemize}

\parSpc
In this paper we present some ideas referring to general tropical
curves and their duals. This includes aspects which subject to the
correspondence that occurs between properties of primal and dual
objects. For a deeper discussion we mainly focus on the family of
quadrics which reflects the brought ideas in a comprehensible
manner and also has significance for other fields of research.
\parSpc
\textbf{\emph{Organization:}} To make this paper reasonably
self-contained we provide a short overview of tropical varieties
and their linkage to complex algebraic varieties. This overview
includes also the construction of the associate Newton polytopes.
Using this basis, in section \ref{sec:PreTropicalDuality}  we
discuss the duality of tropical polynomials (denoted as
pre-tropical duality). Advancing this idea, the
 duality of tropical curves is introduced in section
\ref{sec:TropicalDuality}. We close (Sec.
\ref{sec:RegularityInDualitySense}) by defining the  regularity of
curves in the duality sense and exam it against families of
quadrics. The appendix contains the detailed definitions and
constructions of the classical duality in general case.
Specifically, for the case of quadrics the dual map in matrices'
notion is developed (Sec.
 \ref{subSec:DualityOfMultidimensionalQuadrics}).

% To make this paper reasonably self-contained we provide a short
% overview of ....
% For fluent reading, we open with a review on the fundamentals of
% Non-Archimedean Amoeba.

%************************* Acknowledgement ****************************
\parSpc
\textbf{\emph{Acknowledgement:}} The author would like to thank
Prof. Eugenii Shustin\footnote{The Department of Pure Mathematics,
The School of Mathematical Sciences, Tel Aviv University.} for his
invaluable help. I'm deeply grateful to him for his support.
%\end{AMS}

%******************************* SECT  *********************************
\secSpc
\section{Tropical Varieties}
%\section{Non-Archimedean Amoebaes and Newton Polytope}
For fluent reading, we open with a review of the fundamentals of
Tropical Varieties and Non-Archimedean Amoeba, these fundamentals
are spread among several different works
\cite{Joswig2068,Kapranov2000,Mikhalkin8225,Mikhalkin05,Mikhalkin01,Passare00,Shustin1278}.
%~~~~~~~~~~~~~~~~~~~~~~~~~~~ SUB-SECT ~~~~~~~~~~~~~~~~~~~~~~~~~~~~~~~~~~~
\subSecSpc
\subsection{General Tropical Variety (or Non-Archimedean Amoeba)}
Let $\Fld$ be an algebraically closed field with a valuation
$$\Val : \; \Fld \; \To \; \maxPlus. $$
For our concern the field $\Fld$ is assumed to be the field of
convergent \emph{Puiseux Series}, over the complex numbers
$\Comp$, of the form
\begin{equation}\label{eq:powerSeriesBound}
  a(t) = \sum_{\tau \in R} c_{\tau} t ^{\tau},
\end{equation}
where $R \subset \Rati$ is bounded from below and the elements of
$R$ have a bounded denominator. More specifically, $R$ is
contained in the sum of finitely many arithmetic progressions
which are bounded from below and satisfy,
$$ \sum_{\tau \in R} |c_{\tau}| t ^{\tau} < \infty,$$
for sufficiently small positive $t \neq 0$.

\parSpc
The induced \emph{Non-Archimedean Valuation} over this field is
defined as minus the smallest $\tau$ for which the coefficient
$c_{\tau} \neq 0$, formally,
\begin{equation}\label{eq:valPowerSeries}
  \Val(a(t)) = - min \{\tau \in R \ : \; c_{\tau} \neq 0 \}.
\end{equation}
This valuation takes $\Fld^*$ onto $\Rati$ and one can verify that
for any $a,b \in \Fld^*$ it satisfies the required relations,
\begin{equation}\label{eq:valRules}
\begin{array}{lll}
  \Val(a \cdot b) & = & \Val(a) + \Val(b), \\  [2mm]
  \Val(a+b) & \leq & max \{ \Val(a), \Val(b)\},
\end{array}
\end{equation}
of being Non-Archimedean. Thus, the valuation maps a series $a(t)$
into an element of the semi-ring which defined over $\maxPlus$.
Further on, this semi-ring will be described explicitly. Note
that, the symbol $\tUniS$ is used as the image of the zero element
of $\Fld$ under the valuation \Ref{eq:valPowerSeries}.

\parSpc
Based on the above valuation \label{eq:valRules} we define of the
\emph{Tropical Varieties} (i.e. the \emph{Non-Archimedean Amoeba})
as the closure in $\Real^n$ of the image of some ``superior''
algebraic varieties that were placed on $\Fld^n$. This image
corresponds to the map defined by $\Val$ while the primal elements
(i.e. the algebraic variety) are contained in the zero set of a
system of \emph{Laurent Polynomials} with coefficients in
$\Fld^*$. Formally, let $I \subset \Int^n$ be a non-empty set,
denote by $F_{\Fld}(I)$ the family of Laurent polynomials with
coefficients in $\Fld^*$ having the  form,
\begin{equation}\label{eq:LaurentPolynomials}
  f(\nXz) = \sum_{\om \in I} c_{\om}  \nXz^{\om},
\end{equation}
where $\nXz$ stands for the $n$-tuple $(z_1, \dots, z_n)$. Let,
$f$ be a polynomial in  $F_{\Fld}(I)$, and let
$$Z_f = \{ \nXz \; | \; f(\nXz) = 0 \} \subset (\Fld^*)^n,$$
be its zero locus. The Non-Archimedean amoeba $\Amb_f$ that
corresponds to $f$ is defined to be
\begin{equation}\label{eq:amoebaDef}
 \Amb_f = \overline{\Val(Z_f)} \subset \Real^n,
\end{equation}
where $\Val(z_1,\dots,z_n) = (\Val(z_1),\dots,\Val(z_n))$ and,
$\overline{\Val(Z_f)}$ is the closure of $\Val(Z_f)$. The set of
all amoebas $\Amb_f$, while $f \in F_{\Fld}(I)$, is denoted by
$\Amb(I)$. In case $I$ is the complete set of the integral points
included in a given lattice polygon $\Polygon$ the amoebas' set is
signed by $\Amb(\Polygon)$.

%~~~~~~~~~~~~~~~~~~~~~~~~~~~ SUB-SECT ~~~~~~~~~~~~~~~~~~~~~~~~~~~~~~~~~~~
\subSecSpc
\subsection{Tropical Hyper-Surfaces}

Regarding tropical hyper-surfaces, the first fundamental result
that is brought is due to Kapranov \cite{Kapranov2000}. This
result associates between ``classical'' algebraic objects and
``tropical'' ones (the proof can also be found in
\cite{Einsiedler8311} or \cite{Shustin1278}).
\begin{theorem}\label{thm:Kapranov} {\bf(Kapranov)}. The Non-Archimedean amoeba $\Amb_f$
coincides with the corner locus of the convex piece-wise affine
linear function
$$\nAmb_f (\nXx) = \max_{\om \in I}(\Val(c_\om) + \om.\nXx ),$$
where $\nXx \in \Real^n$.
\end{theorem}
\noindent \textbf{\emph{Notation \& terminology:}}
\begin{itemize}
  \item The term corner locus means domain of non-smoothness.
  \item Unless otherwise specified, the product of the vectors as
  it
  appears here (i.e. $\om.\nXx$),  is holds as the standard scalar product.
\end{itemize}
Moreover, Kapranov's  theorem implies that not only that every
tropical hyper-surface is a corner locus of a tropical polynomial,
but also that any corner locus of a tropical polynomial $\tF$ is a
tropical hyper-surface: simple take $f = \sum_{\om \in I} c_\om
\nXz ^ \om$ such that $\Val(c_\om ) = a_\om$, this kind of
setting is always possible.

\parSpc
The above theorem lays the foundation for understanding the means
of what is called ``\emph{Tropical Geometry}''. By this manner,
$\nAmb$ is regarded as a ``polynomial'' description of a geometric
structure which is embedded in the semi-ring defined over
$\maxPlus$. Hence, $\Real$ can be considered as associated with
the two following ``tropical operations'',
$$
 x \TrS y = \max\{x, y\}, \qquad x \TrP y = x + y.
$$
The triple $\maxPlusSR$ is referred to as the tropical semi-ring.

\parSpc
\begin{remark}
In the above definition of the tropical addition, one may
equivalently replace the $max$ operation by $min$ operation.
However, as will be seen later on, using the operation $max$
appeares to be more convenient for the duality considerations,
especially when the corresponding geometric objects are
constructed in respect to Newton Polytope.
\end{remark}

\begin{figure}[!h]
\centering
\includegraphics[width=\FigWidth in]{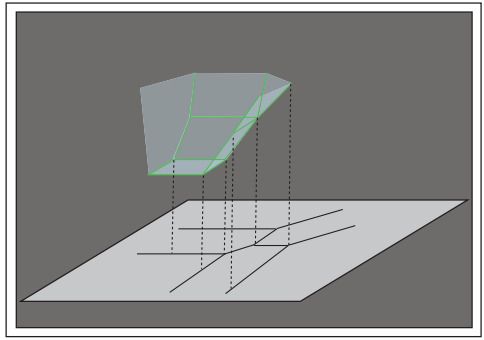}
\caption{\label{fig:cornerLocus} The tropical curve is the
projection of the corner locus of $\Polygon_{\nuf}$, which is
``placed'' above.}
\end{figure}

\parSpc
Composing all the above settings together, the obtained result is
a \emph{Tropical Polynomials} in $n$ variables which are defined
over $\maxPlusSR$, and have  the form
\begin{equation}\label{eq:TropPolynomials}
\tF(\nXx) = \bigoplus_{\om \in I} a_\om \TrP \nXx^ \om,
\end{equation}
where $x_i^{\om_i}$ stands for $x_i \TrP \dots \TrP x_i$, $\om_i$
times. Note that a tropical polynomial can actually be interpreted
(using the ``regular'' operation of addition and multiplication)
as
$$
\tF(\nXx) = \max_{\om \in I}( a_\om +  \nXx . \om),
$$
which is a convex function on $\Real^n$. As a result $\tF$ can be
considered as functional $\Real^n \rightarrow \Real$ and one can
observe its corner locus (see Fig. \ref{fig:cornerLocus}).

\begin{comment}
Geometrically, using the above notion, polynomial descriptions of
piece-wise linear objects are produced. These descriptions will
serve us later in order to manipulate the corresponding objects.
\end{comment}
%%%%%%%%%%%%%%%%%%%%%%%%%%%%%%%%%%%%%%%%%%%%%%%%%%%%%%%%%%%%%%%%%%%%
%%
%%%%%%%%%%%%%%%%%%%%%%%%%%%%%%%%%%%%%%%%%%%%%%%%%%%%%%%%%%%%%%%%%%%%

%~~~~~~~~~~~~~~~~~~~~~~~~~~~ SUB-SECT ~~~~~~~~~~~~~~~~~~~~~~~~~~~~~~~~~~~
\subSecSpc
\subsection{Amoebas, Newton Polytope and Subdivisions}
For a given polynomial equation (either, tropical or classical)
over a set of indices $I \subset \Int^n$ the Newton Polytope
$\Polygon$ (polygon in the planar case) is defined to be the
convex hull, $\CH(I)$, of all the elements of $I$
\cite{Gelfand94}.  Note that, when concerting a tropical
polynomial in respect to the previous definition, all the
coefficients take part. Particularly, zero coefficients should
also appear. As will be explained, this construction of
\emph{Newton Polytopes} validates the linkage between the
algebraic entities and their corresponding geometry objects.

\parSpc
Let $f \in F_\Fld (I)$ be any polynomial in $n$ variable of the
form \Ref{eq:TropPolynomials}, and let $\Polygon_f \subset
\Real^n$ be its corresponding Newton polytope, then $f$ induces a
natural subdivision $\tS_f$ on $\Polygon_f$. This subdivision is
obtained as follows. First, over the points of $\Polygon_f$ define
the \emph{Lifting Map}
\begin{equation}\label{eq:liftingMap}
   \nuf: \Polygon_f \To \Real,
\end{equation}
where, $\nuf(\om)=  -\Val(c_\om) = -a_\om$ (remember that the
$\om$'s are integral points). Then, construct the set
\begin{equation}\label{eq:liftingBody}
J_f = \{ (\om,u) \in I \times \Real \; | \; u > \nuf(\om) \},
\end{equation}
this set is contained in $\Real^{n+1}$. Next, compute the convex
hull, $\Polygon_{\nuf}$, of $J$ and define the function
\begin{equation}\label{eq:convexMap}
   \tilNuf: \Polygon_f \To \Real,
\end{equation}
according to the relation
$$ \tilNuf(\om) = \min \{ u \; |\;  (\om,u) \in \Polygon_{\nuf}\}.$$
This is a convex piece-wise linear function whose domain is the
preliminary polytope $\Polygon_f$ while its vertices are in $I$.
Note that, when an $\om \in I$ is placed on the ``lower part'' of
$\Polygon_{\nuf}$ we have $\tilNuf(\om) = \nuf(\om) = -a_\om$.

\begin{figure}[!h]
\centering
\includegraphics[width=\FigWidth in]{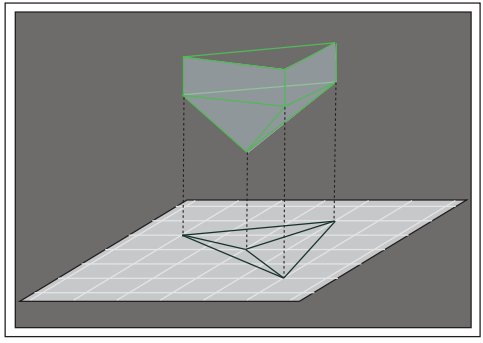}
\caption{\label{fig:polySubdivision} The subdivision induced by a
tropical polynomial.}
\end{figure}

\parSpc Eventually, in order to obtain the required subdivision $\tS_f$
of $\Polygon_f$, take the ``lower part'' $\lbCH(\Polygon_{\nuf})$
of this convex hull, i.e. the graph of the function $\tilNuf$, and
project its non linear part onto $\Polygon_f$ to obtain $\tS_f$ --
the \emph{Polytope's Subdivision} is induced by $f$ (see Fig.
\ref{fig:polySubdivision}).

%In the reminder, where is no concern to confusion we may use the
%simpler notation $\Polygon_{f}$ for $\Polygon_{\nuf}$.

%%%%%%%%%%%%%%%%%%%%%%%%%%%%%%%%%%%%%%%%%%%%%%%%%%%%%%%%%%%%%%%%%%%%
%%
%%%%%%%%%%%%%%%%%%%%%%%%%%%%%%%%%%%%%%%%%%%%%%%%%%%%%%%%%%%%%%%%%%%%

\parSpc
Using the previous construction of the subdivision we can link
between the algebraic description and the geometric notion of the
corresponding tropical object. Globally, there are three different
kinds of entities with reciprocal relationships which play a role
in this view. The ``types'' of entities are, polynomial equations,
polyhedral objects and tropical curves. The linkage between the
two first ``types'' has been already described, and it remains to
understand how do the tropical objects fit in. For this part we
will be assisted by subdivisions of Newton polytopes which have
been defined previously.

\parSpc
In general, any component of a subdivided polyhedron is a
polyhedron by itself. Let $\PolyHedA$ be such a component in
respect to the subdivision $\tS_f$ of $\Polygon_f$ (i.e.
$\PolyHedA \subseteq \Polygon_f$) and let $\PolyHedA_\nu \subseteq
\Polygon_{\nuf} $ be the polyhedron above it (i.e.
$Proj(\PolyHedA_\nu) = \PolyHedA$). Assume  $\PolyHedA_\nu$ and
$\PolyHedB_\nu$ to be two polyhedrons of the same subdivision
$\tS_f$ such that $\PolyHedA_\nu \subset \PolyHedB_\nu$, and
denote by
$$
Cone(\PolyHedB_\nu, \PolyHedA_\nu) = \{t(\nXx -  \nXy) \; | \; t
\geq 0, x \in \PolyHedA_\nu, y \in \PolyHedB_\nu \}, $$
the cone of vectors pointing from $\PolyHedB_\nu$ into
$\PolyHedA_\nu$. Using regular duality in cones' sense (for more
details see \cite{fulton79,oda85}), the dual set can be defined as
$$ DualCone(\PolyHedB_\nu, \PolyHedA_\nu) = \{ \xi \in \Real^{n+1} \; |
  \; \xi.\nXx \leq 0, \;  \forall \nXx \in  Cone(\PolyHedB_\nu, \PolyHedA_\nu) \}.$$
This set appeared to be also a proper cone. Clearly, using the
same manner, such dual cones can be defined for any pair of
polyhedrons compounded in the subdivision $\tS_f$.

\parSpc
Since $\Polygon_{\nuf}$ is by itself a cone (i.e the maximal
cone), we can apply the above considerations to obtain its common
dual cone, $DualCone(\PolyHedA_\nu, \Polygon_{\nuf})$, with
respect to any ``smaller'' $\PolyHedA_\nu$. One can easily
verified that this cone is contained in the set $\{ x_{n+1} \leq 0
\}$ and intersects the hyper-plane $\{ x_{n+1} = -1 \}$. The
$Cone(\PolyHedA_\nu, \Polygon_{\nuf})$ is contained in a half
space (the cone of a non vertical $n$-face of $\Polygon_{\nuf}$
containing $\PolyHedA_\nu$), which is always the half-space above
the (non-vertical) $n$-space cutting $\Real^{n+1}$. Thus,
$DualCone(\PolyHedA_\nu, \Polygon_{\nuf})$ contains at least the
dual of this half space, which is the semi-line intersecting
$\{x_{n+1} = -1\}$.

%%%%%%%%%%%%%%%%%%%%%%%%%%%%%%%%%%%%%%%%%%%%%%%%%%%%%%%%%%%%%%%%%%%%
%%
\begin{figure}[!h]
\centering
\includegraphics[width=\FigWidth in]{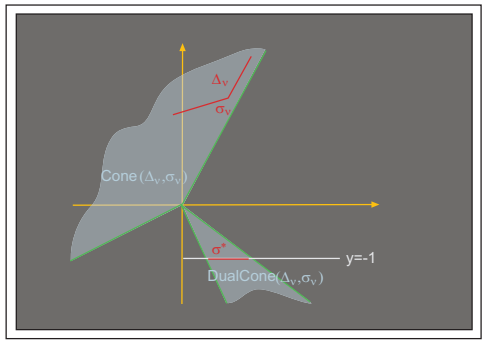}
\caption{\label{fig:dualCones} A cone and its dual.}
\end{figure}
%%%%%%%%%%%%%%%%%%%%%%%%%%%%%%%%%%%%%%%%%%%%%%%%%%%%%%%%%%%%%%%%%%%%

\parSpc
Additionally, from the definition of $\Polygon_{\nuf}$,
$Cone(\PolyHedA_\nu, \Polygon_{\nuf})$ always contains a vertical
semi-line (upper part), hence, $DualCone(\PolyHedA_\nu,
\Polygon_{\nuf})$ is comprised in the half space $\{x_{n+1} \leq
0\}$. Thus a corresponding polyhedron
$$
\PolyHedA^* = \{ \nXx \in \Real^n \; | \; (\nXx,-1) \in
DualCone(\PolyHedA_\nu, \Polygon_{\nuf}) \}
$$
can be defined. Fig. \ref{fig:dualCones} shows this relation for
one dimensional space. In fact, with this construction, a
subdivision of $\Real^n$ is defined, which in the cones' sense is
dual to the subdivision $S_f$ of $\Polygon_f$. Namely, $\PolyHedA
\subset \PolyHedB$ if and only if $\PolyHedB^* \subset
\PolyHedA^*$ and for this case $Cone(\PolyHedA, \PolyHedB)$ is the
dual of $DualCone(\PolyHedB^*, \PolyHedA^*)$. Here, finally, the
required connection between the different objects is cleared up.
\begin{lemma}\label{lem:liftingLegendre}
The polyhedral subdivision of $\Real^n$ given by the $\PolyHedA^*$
where $\PolyHedA$ ranges through the polyhedra of the subdivision
$\tS_f$ of $\Polygon_f$ is exactly the subdivision obtained  by
the corner locus of the Legendre transform of $\nuf$, which is
$$\tilNuf(\nXx) = \max_{\om \in I} \{\om.\nXx -
\nuf(\om)\}.$$
\end{lemma}

%%%%%%%%%%%%%%%%%%%%%%%%%%%%%%%%%%%%%%%%%%%%%%%%%%%%%%%%%%%%%%%%%%%%
%%
% \ifdef{\printFig}{
% \begin{figure}[!h]
% \centering
% \includegraphics[width=\FigWidth in]{LegendreTrans.eps}
% \caption{\label{fig:LegendreTrans}}
% \end{figure}}{}
%%%%%%%%%%%%%%%%%%%%%%%%%%%%%%%%%%%%%%%%%%%%%%%%%%%%%%%%%%%%%%%%%%%%

\begin{figure}[!h]
\centering
\includegraphics[width=\FigWidth in]{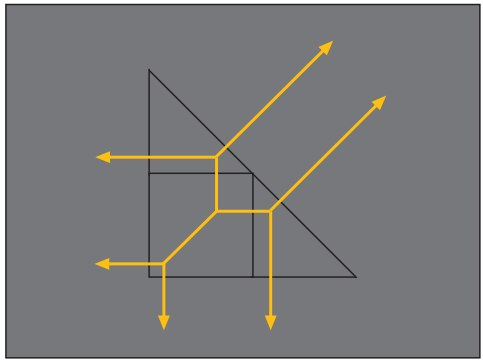}
\caption{\label{fig:amoebaConicNewton} The Newton polytope of
tropical polynomial, its subdivision and, the compatible tropical
conic.}
\end{figure}

\parSpc
It has already been shown that $\tilNuf(\nXx)$  is exactly the
evaluation of the tropical polynomial $\tF$ in $\nXx$. This means
that the subdivision of $\Real^n$ is determined by the tropical
varieties correspond to $\tF$ is dual to the subdivision of the
compatible Newton polytope of $\tF$, see Fig.
\ref{fig:amoebaConicNewton}. One can also see that $\PolyHedA^*$
is orthogonal to $\PolyHedA$ while the dimension of $\PolyHedA$
equals the codimension of $\PolyHedA^*$ (the detailed proof
appears in \cite{Passare00}).

\begin{remark}
Concerning the above geometric observation, characterizing
relations over subdivisions of Newton polytope are equivalent to
specifying relations over tropical objects. Hence, during our
discussion we can base ourself on this insight and develop our
theory in respect to polyhedral objects (i.e. components of
subdivisions of Newton polytope) where the exact translation to
tropical objects is that which is described above.
\end{remark}

%%%%%%%%%%%%%%%%%%%%%%%%%%%%%%%%%%%%%%%%%%%%%%%%%%%%%%%%%%%%%%%%%%%%
%%
%%%%%%%%%%%%%%%%%%%%%%%%%%%%%%%%%%%%%%%%%%%%%%%%%%%%%%%%%%%%%%%%%%%%

%******************************* SECT *********************************
\section{Pre--Tropical Duality}\label{sec:PreTropicalDuality}

As will be discussed later, the tropical duality is practically
induced by the duality which occurs between algebraic verities. In
general, any polynomial $f$ over $\Fld$ has a dual polynomial
$\df$, each has its own lifting (i.e. support) function $\nuf:
\Polygon_f \rightarrow \Real$  (res. $\nuDf: \Polygon_{\df}
\rightarrow \Real$) defined over the integral points of its Newton
polytope $\Polygon_f$ (res. $\Polygon_{\df}$). Using this notion
we establish  a duality between the lifting functions and between
their subdivisions $\tS_f$ ($\tS_\df$) respectively. Moreover,
since the values of these lifting functions over integral points
are the coefficients of the compatible tropical polynomials $\tF$
and $\dtF$, a duality between tropical polynomials is eventually
attained. This duality is called \emph{Pre-Tropical Duality} while
the tropical duality is the duality which occurs between their
corresponding convex functions $\tilde{\nu}_{f}$ and
$\tilde{\nu}_{f^*}$. The last duality can be understood  as the
result of applying convexity considerations to the pre-tropical
duality. At the very end, we will show that the primal algebraic
duality can be disregarded, and we can remain with duality in
tropical terms only.

%~~~~~~~~~~~~~~~~~~~~~~~~~~~ SUB-SECT ~~~~~~~~~~~~~~~~~~~~~~~~~~~~~~~~~~~
\secSpc
\subsection{Pre-Tropical Duality -- Construction}
Let $\tF$ be a tropical polynomial of the form
\Ref{eq:TropPolynomials}, we intend to define its dual $\dtF$. In
order to apply the above considerations, we should specify a
corresponding polynomial $f \in F_{\Fld}(I)$ of the form
\Ref{eq:LaurentPolynomials} such that the coefficients of $\tF$
are the valuations of the coefficients of $f$ (Theorem
\ref{thm:Kapranov}). This, specification can be realized as the
``inverse'' valuation. Recall that, the values of these valuations
are the images of the lifting function $\nuf$ \Ref{eq:liftingMap}
in respect to the points of $I$. However, concerning the
definition of the Non-Archimedean valuation
\Ref{eq:valPowerSeries} it is clear that such $f$ in not unique
and there is a complete family of such polynomials. Moreover, with
this limited information (i.e. the valuations' results) we can't
``recover'' an explicit equation of a polynomial in $\Fld$. But as
will be seen, due to the reapplying of the valuation, the precise
determination of such $f$ is not necessary, and we can avoid it.
The general view over the construction we intent to introduce can
be outlined by the following diagram:
\begin{equation*}
  \begin{CD}
    f                  @>Dual>>         f^*   \\
    @A \iVal AA     @                V {\Val} VV   \\
    \nuf              @> >>            \nu_{f^*}  \\
     @A AA                              @V  VV     \\
     \tF              @> >>             \tF^*
 \end{CD}
\end{equation*}
where the ``missing'' part is $\iVal$.

\parSpc
Before getting on to our development let us first emphasize an
additional important property that relates to Non-Archimedean
valuation.
\begin{remark}\label{rmk:polyValuation} Let
$Q = \sum_{\xi \in J}\pm {\bf c}^\xi$ be a polynomial, by applying
the valuation to $F$ the following relation is obtained
\begin{equation}\label{eq:polyValuation}
  \Val(F) = \Val(\sum_{\xi \in J} \pm {\bf c}^\xi) \leq
\bigoplus_{\xi \in J} \Val( {\bf  c}^\xi).
\end{equation}
In the case that only single component achieves the maximum, the
weak inequality is inverted to be equality.
\end{remark}

\parSpc
Now let us begin our discussion. Let $\tF$ be a tropical
polynomial of the form \Ref{eq:TropPolynomials}, the lifting
function $\nu_\dt: \Polygon_{\nu} \to \Real$ is defined as
$\nu_\dt(\om) = -a_\om$ for any $\om \in I$. Note that, the sign
$``\dt''$ is being used since we don't have an explicit polynomial
$f$. According to the previous values we have to produce a Laurent
polynomial $f \in F_{\Fld}(I)$ of the from
\Ref{eq:LaurentPolynomials} where each of its coefficients,
$c_\om$, satisfies $\Val(c_\om ) = a_\om$. Namely, any ``tropical
coefficient'' $a_\om$ should be mapped as following
$$a_\om \Relate c_\om = \sum_{\tau \in R} c^{(\om)}_\tau t ^{\tau},$$
where $\Val(c_\om ) = a_\om$.

\parSpc
Indeed, in order to achieve the above, an additional requirement
should be imposed, the coefficients $c_\om$'s of $f$ should be
\emph{generic}.  Specifically, generic means that any two
$c_{\om_1}$ and $c_{\om_2}$ which satisfy
$$\Val(c_{\om_1}) = \Val(c_{\om_2}) = - min \{\tau
\in R \ : \; c_{\tau} \neq 0 \}$$
have different coefficients in their leading monomials. Namely,
$c^{(\om_1)}_{\tau} \neq c^{(\om_2)}_{\tau} $ for $\tau = -a_\om$.
Moreover, any combination in $c_{\om_i}$'s dose not zero their
leading monomials. This requirement meant to enable an accurate
transform to the tropical semi-ring even when combined expressions
over these  $c_\om$'s are concerned.

\parSpc The above restriction also insures that the relations which are satisfied by
the valuation \Ref{eq:valRules} are all interpreted as strong
equalities. Obviously, such generic specification of $f$ is
possible and it is not unique. Thus, once we select such generic
polynomial $f$ it is to be regarded as the representative of the
complete family of such polynomials having generic coefficients
$c_\om$'s and satisfy  $\Val(c_\om ) = a_\om$. Therefore, using
such representative $f$, the following correspondences can be
outlined
$$
\tF \; \To \; \nu_\dt \To f.
$$
Namely, $f$ can be related to as the image of $\tF$ under the
``inverse'' valuation $\iVal$ (i.e. $\iVal(\tF) = f$).

\parSpc
In the next step, one can apply the ``classical'' dual
transform\footnote{For those how are not familiar with this topic,
a detailed description, including the explicit construction,
appears in the appendix.} to $f$, i.e. $Dual : f \longrightarrow
\df$. Once, the dual $\df$ had been characterized, the valuation
can be applied to $\df$ which will yield $\dtF$. Let us develop
this part explicitly. Let $\df$ be the dual polynomial of the
form,
\begin{equation}\label{eq:dualLaurentPolynomials}
  \df(\nXz) = \sum_{\dInx \in I^*} c^*_{\dInx}  \nXz^{\dInx},
\end{equation}
by the properties of ``classical'' duality, any $c^*_{\dInx}$ is a
polynomial in $c_{\om}$'s (i.e. $c^{*}_\dInx \in Poly_{\om \in
I}(c_\om)$ where $c_{\om} \in \Fld$). Assume  $m = |I \cap \Int|$
then, by applying the Non-Archimedean valuation to $c^{*}_{\dInx}$
we have,
$$\Val(c^{*}_{\dInx}) = \Val(\sum_{\xi \in J} \pm {\bf c}^\xi) \leq
\bigoplus_{\xi \in J} \Val( {\bf  c}^\xi) = \bigoplus_{\xi \in J}
\bigodot_{\om \in I} \xi_\om \Val(c_\om),$$
where $J \subset \Int^m$ and $ \textbf{c}= (c_{\om_1},\dots,
c_{\om_m})$. Since the preliminary selection of $f$ was to be
generic, referring to Remark \ref{rmk:polyValuation},  the weak
inequality is inverted to equality. Let $\Val(c_\om)=a_\om$, then
the above equation can be written as
\begin{equation}\label{eq:dualTropCoeff}
\Val(c^{*}_{\dInx}) =  \bigoplus_{\xi \in J} \bigodot_{\om \in I}
\xi_\om a_\om = a^*_{\dInx},
\end{equation}
where $a^*_{\dInx}$ is a proper tropical equation in $a_\om$'s for
any $\dInx \in J$. Thus, the valuation of any coefficient
$c^*_{\dInx}$ of $\df$ can be described in tropical terms of the
valuations of the coefficients $c_{\om}$ of $f$. Thus, the
involved values are only the $a_\om$'s up to the multiplications
by integers.
\begin{corollary}
Any
$$a^{*}_{\dInx} \in TropPoly_{\om \in I}(a_\om)$$
and thus we construct a direct map
$$ \mu:\{a_\om \;| \; \om \in I  \} \To \{a^{*}_{\dInx} \;| \; \dInx \in J \}.$$
\end{corollary}
\noindent Clearly, this map is independent in a specific selection
of a representative $f$ unless it has generic coefficients.
Eventually, according to $\mu$ we can define the tropical duality
$$ Dual: \tF \; \longrightarrow \tF^*. $$
where
\begin{equation}\label{eq:dualTropPolynomials}
\dtF(\nXx) = \bigoplus_{\dInx \in I} a^{*}_{\dInx}\TrP \nXx^
\dInx.
\end{equation}
\begin{observation}
In fact, we have shown that in order to characterize the $\dtF$ we
can disregard the returning ``backward'' to computation over the
elements of $\Fld$ and remain only in the tropical semi-ring. In
particular, the explicit construction of $\mu$ is done similarly
to that which is applied to the members of $\Fld$ (for the
complete algorithm, see Appendix) only by replacing the standard
operations by tropical ones.
\end{observation}

Based on the above observation, assume that  $c^{*}_\dInx \in
Poly_{\om \in I}(c_\om)$ was computed, we will outline the rules
of operations' translation,
$$Val: c^{*}_{\dInx} \To a^{*}_{\dInx} $$
where $a^{*}_{\dInx} \in TropPoly_{\om \in I}(a_\om)$. In general,
the basic operations are switched by tropical ones,
$$
\begin{array}{lll}
  c_{\om_1} + c_{\om_2} & \dashrightarrow  & a_{\om_1} \TrS
  a_{\om_2},\\%[2mm]
  c_{\om_1} \cdot c_{\om_2} & \dashrightarrow  & a_{\om_1} \TrP a_{\om_2},
\end{array}
$$
where $\Val(c_{\om_i}) = a_{\om_i}$. Note that any constant
multiplier (including negative signs) of $c_{\om_i}$ is neglected,
and the following maps are valid
$$
\begin{array}{lll}
  n * c_{\om} & \dashrightarrow  & a_{\om},\\ % [2mm]
  -c_{\om}    & \dashrightarrow  & a_{\om},
\end{array}
$$
where $n \in \Int$.
\begin{conclusion}\label{con:TropicalEquivalent}
  Two polynomials in the tropical view are equivalent up to the
  tropical addition of similar monomials, namely any sequence of
  similar monomials,
\begin{equation}\label{eq:topicMomEqu}
 (a_{\om_o} \TrP \nXx^{\om_o}) \TrS \dots \TrS( a_{\om_o} \TrP
 \nXx^{\om_o})
 \simeq a_{\om_o} \TrP \nXx^{\om_o},
\end{equation} can be regarded as a single appearance.
\end{conclusion}
\noindent Note that in this case, since polynomials equations are
involved, this consideration is allowed. This operation is
 interpreted differently while constant values are concerned, and in
that case the tropical addition of equal values is regarded as
minus infinity.

\parSpc
To summarize, using the previous rules, we can perform all the
calculations over the tropical semi-ring. This occurs after
applying the valuation to dual polynomial of that have been
constructed using the ``inverse valuation'',
$$ \Val(\; Poly_{\om \in
I}(\iVal(a_\om))\;) = \Val(\; Daul(\iVal(a_\om))\;).$$
In the next section these ideas are demonstrated for conics.

\parSpc
\begin{example}\label{exmp:preDualPoly} Let $\tF$ be the tropical
polynomial of the form
$$\tF= a_{1} \TrP x_1 \TrP x_2^{2} \TrS a_{2}.$$
Apply $\iVal$ to obtain $f = c_{1}(t)z_1z_2^{2} + c_{2}(t)$ in
$\Fld$ such that $\Val(c_{i}(t)) = a_{i}$, thus $F =
c_{1}(t)z_1z_2^{2}+c_{2}(t)z_0^3$. Substitute
$$\begin{array}{ccc}
  z_0  \mapsto -(\frac{\al_1 u + \al_2 v}{\al_0}),
  & \qquad z_1 \mapsto u, & \qquad z_2 \mapsto v,
\end{array}
$$
so that
$$\begin{array}{ll}
  \tilde{F}(u,v)_{\al_{0},\al_{1},\al_{2}}  = & F(-(\al_1 u + \al_2 v), \al_0 u, \al_0 v) =   \\
  &  c_{2}\al _{1}^{3}u^{3}+3c_{2}\al _{1}^{2}\al _{2}u^{2}v+
  \left( c_{1}\al _{0}^{3}+3c_{2}\al _{1}\al _{2}^{2}\right)
   uv^{2}+c_{2}\al _{2}^{3}v^{3}.
\end{array}
$$
and its two derivatives are
$$
\begin{array}{ll}
  \frac{\partial \tilde{F}}{\partial u} =&
  3c_{2}\al _{1}^{3}u^{2}+6c_{2}\al _{1}^{2}\al _{2}uv+
  \left( c_{1}\al _{0}^{3}+3c_{2}\al _{1}\al _{2}^{2}\right)  v^{2}, \\
  \frac{\partial \tilde{F}}{\partial v} =&
  3c_{2}\al _{1}^{2}\al _{2}u^{2}+2\left( c_{1}\al _{0}^{3}+
   3c_{2}\al _{1}\al _{2}^{2}\right)  vu+3c_{2}\al
   _{2}^{3}v^{2}.
\end{array}
$$
Their resultant is
$$Res(\al_0, \al_1, \al_2)=
3\al _{0}^{6}\al _{1}^{3}\left( 27c_{1}^{2}c_{2}^{2}\al _{1}\al
_{2}^{2}+4c_{1}^{3}c_{2}\al _{0}^{3}\right) $$ neglect the
multiplier and compute
$$\du{f} = Res ( - a_0,  - a_1  , 1) =  27c_{1}^{2}c_{2}^{2}\al _{1}+4c_{1}^{3}c_{2}\al _{0}^{3}.$$
Finally, the dual tropical polynomial $\dtF$ is received by
applying $\Val$ to $c_{i}$,
$$ \dtF = 2a_{1} \TrP 2a_{2} \TrP x_{1} \TrS 3a_{1} \TrP a_{2} \TrP x_{0}^{3}$$
\end{example}
\noindent The analogous computation for the general case
``classical'' can be found in the appendix.
%~~~~~~~~~~~~~~~~~~~~~~~~~~~ SUB-SECT ~~~~~~~~~~~~~~~~~~~~~~~~~~~~~~~~~~~
\subSecSpc
\subsection{Duality of Pre-Tropical Conics}\label{sec:DualityOfPre-TropicalConics}
Next we develop the general form of tropical duality in the case
of conics. Let $\tF$ be a tropical polynomial of degree 2 defined
as
\begin{equation}\label{eq:conicTropMatForm}
 \tF( \bar{x}) = \bar{x} \tA \bar{x}^T =0,
\end{equation}
where $ \bar{x} =  [x_1 ,x_2 , 1 ]$ and $\tA$ is the symmetric
matrix
\begin{equation}\label{eq:conicTropMat}
  \tA = \left[
\begin{array}{ccc}
  a_1 & a_4 & a_5 \\
  a_4 & a_2 & a_6 \\
  a_5 & a_6 & a_3
\end{array} \right].
\end{equation}
%
%Note that the tropical multiplication of vectors is similar to the
%standard one and obey to its rules.
%
%\parSpc
\textbf{\emph{Remark:}} This setting is proper due to conclusion
\ref{con:TropicalEquivalent}, despite the fact that the tropical
monomials $(a_4 \TrP x_1 \TrP x_2 )$ and $(a_5 \TrP x_1 ) \TrS
(a_6 \TrP x_2)$ appear twice.

\parSpc
For this matrix we can fit the matrix
$$  C = \left[
\begin{array}{ccc}
  c_1 & c_4 & c_5 \\
  c_4 & c_2 & c_6 \\
  c_5 & c_6 & c_3
\end{array} \right],
$$
(i.e. $Val(c_i) =a_i$) that corresponds to a polynomial $f$ over
$\Fld$ (see (\ref{eq:conicMat}) in the appendix). Using the same
methods as used for the classical duality, (see appendix
(\ref{eq:conicDualMat})) the dual of $C$ is given as
$$
C^*= \left[ \begin{array}{lll} (c_{6}^{2}-c_{3}c_{2}) & (
c_{6}c_{5} -c_{4}c_{3}) & ( c_{5}c_{2}-c_{4}c_{6}) \\
( c_{6}c_{5} - c_{4}c_{3}) & ( c_{5}^{2}-c_{1}c_{3})  & (
c_{4}c_{5}- c_{1}c_{6}) \\ ( c_{5}c_{2}-c_{4}c_{6}) & ( c_{4}c_{5}
- c_{1}c_{6}) & ( c_{4}^{2}-c_{1}c_{2})
\end{array} \right].
$$
Taking $\Val(C^*)$ with respect to all restrictions mentioned
formerly we get
\begin{equation}\label{eq:conicDualTropMat}
\tA^*= \left[ \begin{array}{lll} (2a_{6} \TrS a_{3}\TrP
a_{2}) & (a_{6} \TrP a_{5}  \TrS a_{4} \TrP a_{3}) & ( a_{5}\TrP a_{2} \TrS a_{4}\TrP a_{6}) \\
( a_{6} \TrP a_{5}  \TrS a_{4} \TrP a_{3}) & ( 2a_{5}  \TrS a_{1}
\TrP a_{3}) & ( a_{4} \TrP a_{5}  \TrS a_{1} \TrP a_{6}) \\ (
a_{5} \TrP a_{2}  \TrS a_{4} \TrP a_{6}) & ( a_{4} \TrP a_{5} \TrS
a_{1}\TrP a_{6}) & ( 2a_{4}  \TrS a_{1} \TrP a_{2})
\end{array} \right].
\end{equation}

\parSpc
This completes the required dual transform and makes the below
diagram commutative:
\begin{equation*}
  \begin{CD}
    C                        @>Dual>>             C^*   \\
      @A \iVal AA        @                    V {\Val} VV   \\
     \tA                        @> >>             \tA^*
 \end{CD}
\end{equation*}
the missing path is constructed by:
$$ \tA \dashrightarrow C \longrightarrow C^* \dashrightarrow \tA^*.$$
and $\tA^*$ defines the dual tropical conics $\dtF$.

\parSpc
In the case of applying the dual transform twice we have an extra
multiplier which is $Det(C)$, this multiplier is now interpreted
in the tropical sense,
\begin{equation}\label{eq:conincsTropDet}
\begin{array}{l}
  \Val(Det(C)) = \tropDet{A} =  \\
  \qquad \max \left(
a_{1}+a_{3}+a_{2},c_{1}+2a_{6},a_{4}+a_{5}+a_{6},a_{2}+2a_{5},2a_{4}+a_{3}\right).
\end{array}
\end{equation}
Note that, the above tropical determinant can be calculated
equivalently by using the primitive tropical operations. This is
done in a similar manner as is done in the classical determinant,
by switching the ``regular'' operations by tropical ones. However,
it might happen that the tropical determinant is invalid, i.e. two
of its components simultaneously achieve the maximum.

%~~~~~~~~~~~~~~~~~~~~~~~~~~~ SUB-SECT ~~~~~~~~~~~~~~~~~~~~~~~~~~~~~~~~~~~
\subSecSpc
\subsection{Pre-Tropical Duality of Quadrics}\label{subSec:Pre-TropicalDualityOfQuadrics}

This part is mostly technical and based on quadrics' duality in
terms of matrices as being developed in appendix
\ref{subSec:DualityOfMultidimensionalQuadrics}. Hence, we will not
get into details, and just review the final result. Let $\tF$ be a
tropical quadric characterized by the tropical matrix $\tA$ then
its dual is defined by $\dtA$ -- the tropical adjoint matrix of
$\tA$.

%******************************* SECT *********************************
\secSpc
\section{Tropical Duality}\label{sec:TropicalDuality}

In general, the tropical duality is obtained by involving
convexity consideration into the pre-tropical duality. As a result
we can discuss the ``behavior'' of dual subdivisions in respect to
their sources, and thus deduce the duality of Non-Archimedean
Amoebas. After we clarify the connection between the classical and
tropical objects for this matter, we can eventually leave the
beyond classical construction and persist only with the tropical
one.

%~~~~~~~~~~~~~~~~~~~~~~~~~~~ SUB-SECT ~~~~~~~~~~~~~~~~~~~~~~~~~~~~~~~~~~~
\subSecSpc
\subsection{Duality of Tropical Objects}

As mention earlier, the pre-tropical duality is used as the  base
for including convexity issues into tropical duality, and it also
gives the linkage to the ``classical'' duality. However, the
convexity issues in this case refer to the properties of the
polyhedra which are defined by the lifting functions and their
corresponding subdivisions. These issues will be studied next.

\parSpc
Assume $\tF$ is a tropical polynomial and $f \in F_{\Fld}(I)$ is
the representative of the polynomials' family whose members have
the valuation correspond to $\tF$. Let us remind ourselves that a
lifting function, $\nuf$, is defined over the Newton polytope
$\Polygon_f$ and the values of $\nuf$ over the integral points
$\om$'s of $\Polygon_f$ relate to the coefficients of $f$, i.e.
$\nuf(\om) = \Val(c_\om)$. But, according to our setting, we have
$\Val(c_\om) = a_\om$ which are the coefficients of the compatible
tropical polynomial $\tF$ and thus $\Polygon_f = \Polygon_\tF$.
According to this insight, we can identify the lifting function,
$\nuf$ with this of $\tF$ and write $\nuf = \nuF$.

\parSpc
Using this notion, the convex function $\tilNuf(\nXx) = \max_{\om
\in I} \{\om.\nXx - \nuf(\om)\}$ can be rephrased in terms of the
coefficients of $\tF$ as
\begin{equation}\label{eq:convexTropFunc}
  \tilNuF(\nXx)= \max_{\om \in I} \{\om.\nXx + a_{\om}\}
\end{equation}
(see Lemma \ref{lem:liftingLegendre}). Thus, the convex
polyhedrons $\Polygon_{\nuf}$ and $\Polygon_{\nuF}$ which are
specified by $\tilNuf(\nXx)$ and $\tilNuF(\nXx)$ are the same one.
Finally, since $\Polygon_{\nuF} =\Polygon_{\nuf}$, both induce the
similar subdivision $S_f = S_\tF$. According the above
``translations'', we actually leave behind any involvement of $f$
and remain only with the tropical polynomial $\tF$.

\parSpc
All the previous ``translations'' have been made for some $f \in
F_{\Fld}(I)$, clearly the same considerations can be applied to
its dual $\df$. The dual relation between $f$ and $\df$ induces
the duality on the entire participants, and particularly the
duality between subdivisions,
\begin{equation}\label{eq:dualSubdivision}
Daul: S_\tF \To S_{\dtF},
\end{equation}
which is derived from the dual relation between $\tilNuF$ and
$\tilDuNuF$. In general this duality is taking pace over
polyhedra, and we have already explained the linkage between these
type of objects to tropical curves. So, in the future we can
concentrate in polyhedra, assuming the map to the tropical curves
is known. However, a global analysis of convexity relations
between polyhedrons is too complicated and we should restrict
ourselves to certain cases. Some of these cases are bought in the
remaining sections.

%~~~~~~~~~~~~~~~~~~~~~~~~~~~ SUB-SECT ~~~~~~~~~~~~~~~~~~~~~~~~~~~~~~~~~~~
\subSecSpc
\subsection{Duality of Tropical Conics}
Based on conics' definition in terms of matrices (i.e. $\tF(
\bar{x}) = \bar{x} \tA \bar{x}^T$, Sec.
\ref{sec:DualityOfPre-TropicalConics}), in order to observe
relationships between any subdivisions  $S_\tF$ and $S_{\dtF}$, we
should first define the map $\tA \To J_\tF$ ($J_\tF$ is defined
according to \Ref{eq:liftingBody}). Specifically, $J_\tF$ can be
realized as $\Polygon_{\tF}$ with attached values over its
integral points, thus for convenience we can regard this map as
having  the correspondence $\tA \To \Polygon_\tF$. The required
correspondence can simply be sketched as following:
\begin{equation}\label{diag:matrixToPolytop}
  \begin{CD}
   \tA  =  \left[
\begin{array}{ccc}
  a_1 & a_4 & a_5 \\
  a_4 & a_2 & a_6 \\
  a_5 & a_6 & a_3
\end{array} \right]                      @> >>
\Polygon_{\nu_\tF} :
\begin{array}{ccclc}
  a_1 &           &     &            & \\
  |   & \diagdown &     &            & \\
  a_5 & -        & a_4 &            & \\
  |   &           & |   & \diagdown  & \\
  a_3 & -        & a_6 & -         & a_2.
\end{array}
 \end{CD}
\end{equation}
Composing this representation with the result we had developed
for the dual conic in matrices form (\ref{eq:conicDualTropMat})
the following is obtained:
%%%%%%%%%%%%% \printExtend %%%%%%%%%%%%%%%
\ifdef{\printExtend}{
\begin{equation}\label{eq:dualPoly2Poly}
  \begin{CD}
\Polygon_{\nu_\tF}:
\begin{array}{ccclc}
  a_1 &           &     &            & \\
  |   & \diagdown &     &            & \\
  a_5 & -        & a_4 &            & \\
  |   &           & |   & \diagdown  & \\
  a_3 & -        & a_6 & -         & a_2
\end{array} \\
                 @V {Dual} VV \\
\Polygon_{\dNu_\tF}^*:
\begin{array}{ccclc}
  2a_{6} \TrS a_{3}a_{2} &           &     &           & \\
  |   & \diagdown &     &            & \\
  a_{5}a_{2} \TrS a_{4}a_{6} & -     & a_{4}a_{3} \TrS a_{6}a_{5} &            & \\
  |   &           & |   & \diagdown  & \\
  2a_{4} \TrS a_{1}a_{2} & -         & a_{1}a_{6} \TrS a_{4}a_{5} & -         & 2a_{5} \TrS a_{1}a_{3}
\end{array}
  \end{CD}
\end{equation}
} {
\begin{equation}\label{eq:dualPoly2Poly}
  \begin{CD}
\begin{array}{ccclc}
  a_1 &           &     &    J_\tF        & \\
  |   & \diagdown &     &            & \\
  a_5 & -        & a_4 &            & \\
  |   &           & |   & \diagdown  & \\
  a_3 & -        & a_6 & -         & a_2
\end{array}
                 @ >{Dual} >>
\begin{array}{ccclc}
  2a_{6} \TrS a_{3}a_{2} &           &     &       J_{\dtF}      & \\
  |   & \diagdown &     &            & \\
  a_{5}a_{2} \TrS a_{4}a_{6} & -     & a_{4}a_{3} \TrS a_{6}a_{5} &            & \\
  |   &           & |   & \diagdown  & \\
  2a_{4} \TrS a_{1}a_{2} & -         & a_{1}a_{6} \TrS a_{4}a_{5} & -         & 2a_{5} \TrS a_{1}a_{3}
\end{array}
  \end{CD}
\end{equation}
}
where the multiplications above are the topical ones. Eventually,
the values that are placed over the integral points determine the
subdivisions  $S_\tF$ and $S_{\dtF}$. These  subdivisions are
determined according to convexity considerations and this issue
will be analyzed in the next section.

\begin{terminology} The diagram of Newton polytopes
with the above form \Ref{eq:dualPoly2Poly} can be referred to as
geometric graphs. Thus, for convenience in the rest of our
discussion, an integral point will be denoted \textbf{Node} and
nodes which are not intermediate points are called
\textbf{Vertices}. In addition, the term \textbf{Edge} refers to a
1-dimensional face (i.e. geometrically a segment) which connects
two nodes.
\end{terminology}
%~~~~~~~~~~~~~~~~~~~~~~~~~~~ SUB-SECT ~~~~~~~~~~~~~~~~~~~~~~~~~~~~~~~~~~~
\subSecSpc
\subsection{Analysis of Conics Duality}
Aiming to analyze convexity properties, which eventually determine
a subdivision, we will need the following property which regards
to arbitrary points whose projections are co-linear. Let $a_i$ and
$a_j$ be two points which are placed on the ``lower part'' of the
convex hull $\lbCH(\Polygon_{\nu_\tF})$. A middle point $a_k \in
\Polygon_{\nu_\tF}$ appears also on $\lbCH(\Polygon_{\nu_\tF})$ if
and only if
$$
  a_k \leq t a_i+(1-t) a_j, \qquad t \in [0,1].
$$
Note that, this property is valid to any arbitrary triple and not
only for internal points.

\parSpc
Now, assume $\tF$ is a conic and apply this observation to points
$a_\om$'s that correspond to integral points of $\Polygon_{\tF}$
that are not vertices. Since the distances between neighboring
integral points of $\Polygon_{\tF}$ are all one, $t$ can be fixed
to be $1/2$, and all the values that correspond to the integral
points can be written in terms of the vertices' values as:
\begin{equation}\label{eq:conicConvexityInTropical}
  \begin{array}{ccc}
  a_4 & = & \frac{1}{2}a_1 \TrP \frac{1}{2}a_2 \TrP \ep_4 \\
  a_5 & = & \frac{1}{2}a_1 \TrP \frac{1}{2}a_3 \TrP \ep_5 \\
  a_6 & = & \frac{1}{2}a_2 \TrP \frac{1}{2}a_3 \TrP \ep_6
\end{array}.
\end{equation}

\begin{observation}\label{obs:epValue}
Using the above notations, any $a_i$ lays on
$\lbCH(\Polygon_{\nuF})$ only when $\ep_i \leq 0$.
\end{observation}

Substituting these values in $\Polygon_{\dNuF}$ according to
\Ref{eq:dualPoly2Poly} yields
$$
\begin{array}{ccclc}
  _{a_2 \TrP a_3  \TrP(2\ep _{6} \TrS 0 )} &           &     &            & \\
  |   & \diagdown &     &            & \\
  _{\frac{1}{2}(a_{1} \TrP 2a_{2} \TrP a_{3})\TrP(\ep _{5}\TrS (\ep _{4} \TrP \ep _{6}))} &
  - &
  _{\frac{1}{2}(a_{1} \TrP a_{2} \TrP 2a_{3})\TrP(\ep _{4}\TrS (\ep _{5} \TrP \ep _{6}))} &            & \\
  |   &           & |   & \diagdown  & \\
  _{a_1 \TrP a_2  \TrP(2\ep _{4} \TrS 0 )} & -         &
  _{\frac{1}{2}(2a_{1} \TrP a_{2} \TrP a_{3})\TrP(\ep _{6}\TrS (\ep _{4} \TrP \ep _{5}))}
  & -         & _{a_1 \TrP a_3  \TrP(2\ep _{5} \TrS 0
  ).}
\end{array}$$
As one can observe, the tropical ``factor'' of any node which is
not a vertex is exactly the average of the ``factors'' of its two
neighboring vertices along a same line. Thus, for the convexity
matter we can disregard these tropical multipliers and retain only
with
\begin{equation}\label{eq:conicsPolyByEp}
  \begin{array}{ccclc}
  \dEp _{1} &           &     &            & \\
  |   & \diagdown &     &            & \\
  \dEp _{5} & - &   \dEp _{4} &            & \\
  |   &           & |   & \diagdown  & \\
 \dEp _{3} & - & \dEp _{6}  & -  & \dEp _{2}
\end{array} \; \; = \; \;
  \begin{array}{ccclc}
  2\ep _{6} \TrS 0  &           &     &            & \\
  |   & \diagdown &     &            & \\
  \ep _{5}\TrS (\ep _{4} \TrP \ep _{6}) &
  - &
  \ep _{4}\TrS (\ep _{5} \TrP \ep _{6}) &            & \\
  |   &           & |   & \diagdown  & \\
  2\ep _{4} \TrS 0 & -         &
  \ep _{6}\TrS (\ep _{4} \TrP \ep _{5})
  & -         & 2\ep _{5} \TrS 0.
\end{array}
\end{equation}
For our purpose this result is called the \emph{``Distortions''}
of $\dNuF$. In addition,  by using
(\ref{eq:conicConvexityInTropical}) this consideration can be
generalized to any $n \times n$ symmetric matrix.

\parSpc
Naturally, the questions that should be asked are, in which
scenarios a point $a^*_{\dInx}$ is contained in
$\lbCH(\Polygon_{\dNuF})$ and, what is the dependency of a dual
coefficient on those who appear in its description. In other
words, what is the induced dual subdivision $S_{\dtF}$ with
respect to the primal subdivision $S_{\tF}$. Clearly this
dependency determines the duality between the induced
subdivisions, $S_{\dtF}$ and its dual $S_{\tF}$. Using this
insight, our discussion can be restricted to subdivisions. For
this matter  we use the notation $S_\Polygon$ to denote a
subdivision of a polytope $\Polygon$.

\parSpc
The immediate observation of \Ref{eq:conicsPolyByEp} yields that,
only nodes of $J_\tF$  which are not vertices involve in
determining whether a node of $J_{\dtF}$ is comprised or not in
$\lbCH(\Polygon_{\dNuF})$. These conditions can be summarized
using the corresponding ``distortions'' values
(\ref{eq:conicsPolyByEp}) as:
$$
\begin{array}{ccc}
  a^*_4 \in \lbCH(\Polygon_{\dNuF}) & \Longleftrightarrow & \ep_4 < (\ep _{5} \TrP \ep _{6}) \TrS \ep _{5} \TrS \ep _{6} \TrS 0,
  \\ [2mm]
  a^*_5 \in \lbCH(\Polygon_{\dNuF}) & \Longleftrightarrow & \ep_5 < (\ep _{4} \TrP \ep _{6}) \TrS \ep _{4} \TrS \ep _{6} \TrS 0,
  \\ [2mm]
  a^*_6 \in \lbCH(\Polygon_{\dNuF}) & \Longleftrightarrow & \ep_6 < (\ep _{4} \TrP \ep _{5}) \TrS \ep _{4} \TrS \ep _{5} \TrS
  0.
  \end{array}
$$
\textbf{\emph{Remark:}} The interesting situations are those when
an image of a node is vertex of  $\lbCH(\Polygon_{\dNuF})$, thus
in the above conditions the inequalities are the strong ones.

\parSpc
\begin{definition}\label{prop:CompleteSubDiv} (\textbf{Complete subdivision})
A subdivision $\tS_\Polygon$ of a polytope $\Polygon$ is called
complete if it can not be refined into other subdivision
$\tS'_\Polygon$.
\end{definition}
\noindent This definition is given in the sense of the resulting
components, namely, none of the subdivision's components can be
subdivided further.  On the other hand the next classification is
defined in nodes sense, and it can be regarded as an ``almost''
complete subdivision.

\parSpc
\begin{definition}\label{prop:FullSubDiv}(\textbf{\full \; subdivision})
A subdivision $\tS_\Polygon$ of a polytope $\Polygon$ is called
\full \; if all the polytope's nodes appear in $\tS_\Polygon$. A
subdivision $\tS_\Polygon$ is called \emty \; if only the vertices
of $\Polygon$ appear in $\tS_\Polygon$.
\end{definition}
\parSpc
Any complete subdivision is always \full \; while the contrary is
not necessarily true. As an example,  chose the following
subdivision,
$$
\begin{array}{ccclc}
  \bullet  &           &     &            & \\
  |  & \diagdown &     &            & \\
  \bullet  &    \ldots     &  \bullet &            & \\
  |  &      \diagdown     & |   & \diagdown  & \\
  \bullet  & -        &  \bullet & -          & \bullet
\end{array}
$$
which is not a ``perfect'' triangulation (i.e. not all of its
components are primitive triangles) and it can be refined, this
despite all of the vertices involved in the subdivision.
\parSpc
\begin{proposition}\label{prop:conicFullSubDiv}
 The dual subdivision $S_{\dtF}$ of
\begin{enumerate}
  \item a subdivision $\tS_\tF$ which is \full \; is \full.
  \item a subdivision $\tS_\tF$ which is \emty \; is \emty.
\end{enumerate}
\end{proposition}
\begin{proof}
The two assertions derived form construction
(\ref{eq:conicsPolyByEp}) using the corresponding $\ep$'s values
as defined in (\ref{eq:conicConvexityInTropical}).
\begin{enumerate}
%%%--------------------------
  \item
  According to setting (\ref{eq:conicConvexityInTropical}) a subdivision
  $\tS_\tF$ is \full \; if and only if
  $ \ep _{i} < 0 $ for all $i=4,5,6$. In respect to
  these values, $\ep _{i} \TrS 0  = 0$, and any combination which
  satisfies the relation  $\dEp _{i} = \ep _{i}\TrS (\ep _{j} \TrP \ep
  _{k}) < 0$ where $i,j,k \in \{4,5,6\}$. Thus, the corresponding dual
  image as described in (\ref{eq:conicsPolyByEp}) can be written as
  \begin{equation*}
  \begin{array}{ccclc}
  0  &           &     &            & \\
  |   & \diagdown &     &            & \\
  \ep _{5}\TrS (\ep _{4} \TrP \ep _{6}) &
  - &
  \ep _{4}\TrS (\ep _{5} \TrP \ep _{6}) &            & \\
  |   &           & |   & \diagdown  & \\
  0 & -         &
  \ep _{6}\TrS (\ep _{4} \TrP \ep _{5})
  & -         &  0.
\end{array} \; \; = \; \;   \begin{array}{ccclc}
  0  &           &     &            & \\
  |   & \diagdown &     &            & \\
  \dEp _{5} < 0& - &  \dEp _{4} < 0&            & \\
  |   &           & |   & \diagdown  & \\
  0 & -  &  \dEp _{6} < 0
  & -         &  0.
\end{array}
\end{equation*}
As one can easily observe, any node which is not vertex has a
negative value and hence, it appears in the lower part
$\lbCH(\Polygon_{\dNuF})$. Since, this is valid  for any node, the
compatible subdivision, $S_{\dtF}$ is \full.
%%%--------------------------
  \item
 As defined in (\ref{eq:conicConvexityInTropical}) a subdivision $\tS_\tF$
 is \emty \; if and only if
  $ \ep _{i} > 0 $ for $i=4,5,6$. In respect to these
  values $\ep _{i} \TrS 0  = \ep _{i}$ and any
  combination satisfies the relation
  $\ep _{i}\TrS (\ep _{j} \TrP \ep
  _{k}) \geq \ep _{j} \TrP \ep
  _{k}$ where $i,j,k \in \{4,5,6\}$. The relative values of the
  corresponding dual polytope can be written via
  (\ref{eq:conicsPolyByEp}) as
  \begin{equation*}
  \begin{array}{ccclc}
  2\ep _{6}  &           &     &            & \\
  |   & \diagdown &     &            & \\
  \ep _{5}\TrS (\ep _{4} \TrP \ep _{6}) &
  - &
  \ep _{4}\TrS (\ep _{5} \TrP \ep _{6}) &            & \\
  |   &           & |   & \diagdown  & \\
  2\ep _{4}  & -         &
  \ep _{6}\TrS (\ep _{4} \TrP \ep _{5})
  & -         & 2\ep _{5}.
\end{array}
\end{equation*}
The obtained result is that, only the vertices appear in the
subdivision which means that it is \emty.
\end{enumerate}
\end{proof}

\parSpc
The two following examples present these kind of subdivisions, and
in addition not only that the dual is shown  but also the result
of applying the duality twice. This and the calculations of the
tropical determinants (\ref{eq:conincsTropDet}) will assist us
 in our later development.
\begin{example}\label{examp:fullSubD}
A subdivision which is \full \; is mapped via the dual tropical
transform (\ref{eq:dualPoly2Poly}) into a subdivision which is
also  \full. In addition the diagram on the left is the result of
applying the duality twice.

$$
\begin{array}{ccclc}
     &       &    \nuF &            & \\
  1  &           &     &            & \\
  |  & \diagdown &     &            & \\
  1  &  -        &  (-2) &            & \\
  |  &           & |   & \diagdown  & \\
  3  & -         &  (-1) & -          & 2
\end{array}
\longrightarrow
\begin{array}{ccclc}
     &    &  \dNuF    &            & \\
  5  &           &     &            & \\
  |  & \diagdown &     &            & \\
  3  &  -        &   1 &            & \\
  |  &           & |   & \diagdown  & \\
  3  & -         &  0  & -          & 4
\end{array}
\longrightarrow
\begin{array}{ccclc}
     &    & \nu_{\ddtF}&            & \\
  7  &           &     &            & \\
  |  & \diagdown &     &            & \\
  7  &  -        &  4  &            & \\
  |  &           & |   & \diagdown  & \\
  9  & -         &  5 & -           & 8.
\end{array}
$$
For this case the tropical multiplier for passing from $\nuF$ to
$\nu_{\ddtF}$  is,
$$TropDet(\tA) =
%\max \left\{
%a_{1}+a_{3}+a_{2},a_{1}+2a_{6},a_{4}+a_{5}+a_{6},a_{2}+2a_{5},2a_{4}+a_{3}\right\}$$
%$$
\max \{1+3+2,1-2,-2+1-1,2+2,-4+3\}= 6$$ and the nodes' values
 are equal up to the lifting by a constant. Thus, $\nuF$ and
$\nu_{\ddtF}$, induce the same subdivision (i.e. $\tS_\tF$ =
$\tS_{\ddtF}$).
\end{example}

\begin{example}\label{examp:emptySubD}
As done in the previous example the three diagrams below refer to
$\nuF$, $\dNuF$ and $\nu_{\ddtF}$ for the case that
 $\nuF$ induces a subdivision which is \emty,
$$
\begin{array}{ccclc}
     &       &    \nuF &            & \\
  1  &           &     &            & \\
  |  & \diagdown &     &            & \\
  4  &  -        &  3 &            & \\
  |  &           & |   & \diagdown  & \\
  5  & -         & 4 & -          & 2
\end{array}
\longrightarrow
\begin{array}{ccclc}
     &    &  \dNuF    &            & \\
  8  &           &     &            & \\
  |  & \diagdown &     &            & \\
  7  &  -        &  8 &            & \\
  |  &           & |   & \diagdown  & \\
  6  & -         &  7  & -          & 8
\end{array}
\longrightarrow
\begin{array}{ccclc}
     &    & \nu_{\ddtF}     &            & \\
  14  &           &     &            & \\
  |  & \diagdown &     &            & \\
  15  &  -        &  14  &            & \\
  |  &           & |   & \diagdown  & \\
  16  & -         &  15 & -           & 14.
\end{array}
$$
For this case the tropical determinant is,
$$TropDet(\tA) =
% \max \left\{
%a_{1}+a_{3}+a_{2},a_{1}+2a_{6},a_{4}+a_{5}+a_{6},a_{2}+2a_{5},2c_{4}+c_{3}\right\}
%$$
%$$
\max \{1+5+2,1+2*4,4+4+3,2+2*4,2*3+5\}= -\infty.$$ However, as can
be seen this does not harm the result, and the dual of subdivision
which is \emty \; is also \emty.
\end{example}

%~~~~~~~~~~~~~~~~~~~~~~~~~~~ SUB-SECT ~~~~~~~~~~~~~~~~~~~~~~~~~~~~~~~~~~~
\subSecSpc
\subsection{Duality of Tropical Quadrics}
In this section we will generalize the previous results of
subdivisions corresponding to conics for the case of subdivisions
that refer to quadrics. As before, we base ourselves on the
linkage between the two classes of associated objects: the
classical and the tropical where both are given in matrices form.
For this purpose, we will develop assisting methods which lead the
further analysis.

\parSpc
As a reminder, let
\begin{equation*}\label{eq:generalQuadrics}
  f(\brZ) = \brZ C \brZ^T=0
\end{equation*}
be a quadric in $\Fld^n$, where $C$ is the $(n+1)\times(n+1)$
symmetric matrix and $\brZ = (z_1, \dots, z_n , 1)$. Its dual is
defined by,
\begin{equation*}\label{eq:generalQuadricsDual}
  \df(\brZ) = \brZ^*\left[C^{-1}\right] (\brZ^*)^T
  =0.
\end{equation*}
For subdivision issues, which are influenced only by convexity
properties, one can rely on the relation
\begin{equation*}\label{eq:conicByInverseMat}
 C^{-1} = Adj(C)\left(\frac{1}{Det(C)}\right),
\end{equation*}
and write
\begin{equation*}\label{eq:generalQuadricsDual}
  f^*(\brZ) = \brZ^*\left[C^*\right] (\brZ^*)^T
  =0
\end{equation*} where $C^*$ denotes the $Adj(C)$ the adjoint matrix.

\parSpc
Recall that, an entry $c_{i,j}^*$ of $C^*$ is equal to the
determinant $Det(C_{i,j})$ of the $i,j$-minor $C_{i,j}$. In
general we have the equation,
\begin{equation}\label{eq:deteminantByMinors}
  Det(C) = \sum^{n+1}_{k=1} (-1)^{h+k} c_{h,k} Det(C_{h,k}).
\end{equation} For convenience $(-1)^{h+k}$ is sometimes signed as
$\pm$ and $Det(C_{h,k}) = |C_{h,k}|$.

\parSpc
Let $V(\tF)$ be a tropical quadric defined by the matrix $\tA$ and
let $f$ be a quadric over $\Fld^n$ defined by the matrix $C$,
where the entries $a_{i,j}$ of $\tA$ are the valuation of
$c_{i,j}$ -- the entries of $C$.  Thus, according to this
construction, $\tA$  can be written as $\tA = \Val(C)$ and
$\Val(Det(C)) = \tropDet{A} = |\tA|$. Clearly, the same holds for
any two minors $\tA_{h,k}$ and $C_{h,k}$. These results coincide
with the results which are obtained by the ``pure'' calculation of
tropical determinants. The basics of tropical determinants and
their properties can be found
\cite{Develin8254,Joswig2068,Sturmfels6366}. Let $\dtA$ be the
matrix with entries $a^*_{i,j} = \Val(c^*_{i,j})$ thus, its
entries are the tropical determinants $|\tA_{i,j}|$. Hence,
eventually, for our concern we can remain with only the two
tropical matrices $\tA$ and $\dtA$, and disregard the beyond
constructions.

\parSpc
In order to observe convexity properties, as done for conics (see
(\ref{diag:matrixToPolytop})), we should first map the matrix
$\tA$ into the corresponding Newton Polytope $\Polygon_\tF$, this
is done using the following rules:
$$
\begin{array}{llll}
  a_{i,i} & \Relate & Zero \; Vertex \; of \; \Polygon_\tF   & i = n+1, \\ [2mm]
  a_{i,i} & \Relate & Vertex \; of \; \Polygon_\tF & i \neq
  n+1, \\[2mm]
  a_{i,j} & \Relate & Node \; of \Polygon_\tF  & i \neq j .
\end{array}
$$
The scheme of this map in respect to one of the polytope's faces
is as follow:
\begin{equation}\label{diag:genMatrixToPolytop}
  \begin{CD}
 \tA:   \left[
\begin{array}{ccc}
  a_{1,1} &  &  \\
  \vdots &  \ddots &  \\
  a_{n,1} & \ldots & a_{n,n}
\end{array} \right]                  @> >>
\Polygon_\tF|_{i,j,k}:
\begin{array}{ccclc}
   a_{j,j}  &           &     &            & \\
  |  & \diagdown &     &            & \\
   a_{j,k}  &    -     &  a_{i,j} &            & \\
  |  &           & |   & \diagdown  & \\
   a_{k,k}  & -        &  a_{i,k} & -          & a_{i,i}.
\end{array}
\end{CD}
\end{equation}

\parSpc
Recall that, in this specific case a node is placed on the middle
of some boundary edge. Thus, the convexity of the lifting function
over these nodes can be determined via the interior relationships
over all the triples,
\begin{equation}\label{eq:tripelsRelatiov}
  \da _{i,i} \longleftrightarrow \da _{i,j} \longleftrightarrow \da _{j,j}.
\end{equation}
For future development the following notation will be needed, let
$[\tA_{i,j}]_{h,k}$ be the $h,k$-minor of the minor $\tA_{i,j}$
such that $ h \neq i$ and $ k \neq j$ with respect to the primal
indices of $\tA$, then $[\tA_{i,j}]_{h,k}$ is signed shortly as
$[\tA_{ih,jk}]$. In the reminder of our deployment, unless
otherwise specified, in order to avoid confusion all indices will
be taken relatively to the primal matrix. In addition, by the
symmetry of $\tA$, the following are also valid:
$$
  \tA_{ih,jk}  =  \tA_{jk,ih}^T, \qquad
  \tA_{ih,jk}  =  \tA_{hi,kj}.
$$
%$$ \begin{array}{ccc}
%  \tA_{ih,jk} & = & \tA_{jk,ih}^T, \\
%  \tA_{ih,jk} & = & \tA_{hi,kj}.
%\end{array}
%$$
%
\parSpc
Using the tropical analogue  to the determinant construction  by
minors over $\Fld$, Eq. (\ref{eq:deteminantByMinors}), and the
above notation any entry $a^*_{i,j}$ can be written as,
\begin{equation}\label{diag:tropiclaConvexity}
  a^*_{i,j} =  \bigoplus_{k\neq j} a_{h,k} \TrP |\tA_{ih,jk}|, \qquad  h \neq i
\end{equation}
%
%\begin{equation}\label{diag:tropiclaConvexity}
%\begin{array}{clll}
%  a^*_{i,i}: & & \bigoplus_{k\neq i} a_{h,k} \TrP |\tA_{ih,ik}|, & h \neq i \\
%  | & & & \\
%  a^*_{i,j}: & & \bigoplus_{k\neq j} a_{h,k} \TrP |\tA_{ih,jk}|, &  h \neq i\\
%  | & &  & \\
%  a^*_{j,j}: &  & \bigoplus_{k\neq j}a_{h,k} \TrP |\tA_{jk,jk}|, &  h \neq j
%\end{array}
%\end{equation}
where $k= 1,\dots,n+1$.

\parSpc
Since, the distances between neighboring nodes of the Newton
Polytope are fixed to one (as was in
(\ref{eq:conicConvexityInTropical})) we generally define the
``distortion'' value $\ep_{i,j}$ of $a_{i,j}$ using the relation
\begin{equation}\label{eq:conicConvexityInTropicalMultiD}
  a_{i,j} =
  \hlf a_{i,i}\TrP \hlf a_{j,j} \TrP  \ep_{i,j},
\end{equation}
where $\ep_{i,i}$ are always zero.

\begin{remark}
The above construction is invariant in respect to convexity, i.e.
a triple is convex if and only if the triple of its
``distortions'' is convex.
\end{remark}

\begin{definition}\label{def:disMat}
Let $\tA$ be a symmetric tropical matrix, the ``distortion''
matrix $\epMat$ of $\tA$ is defined to be $(\ep_{i,j})$ where the
entries $\ep_{i,j}$ are the ``distortion'' values of each entry
which is obtained via (\ref{eq:conicConvexityInTropicalMultiD}).
\end{definition}

\parSpc
Using the ``distortion'' values of $a_{i,j}$, Eq.
(\ref{eq:conicConvexityInTropicalMultiD}), we can represent the
entries $a^*_{i,j}$ of $\dtA$, this will make the convexity
analysis easier. Let $g_{i,j}(\tA)$ be the ``secondary''
expression defined as,
\begin{equation}\label{eq:gFuncTropicalMultiD}
   g_{i,j}(\tA) = \left\{ \begin{array}{lll}
     0, &  & i = j,  \\ [2mm]
     \frac{1}{2}a_{i,i} \TrP \frac{1}{2}a_{j,j},&  & i \neq j.
   \end{array}
   \right.
\end{equation}
and let the ``main'' expression be
\begin{equation}\label{eq:GijFuncTropicalMultiD}
   G_{i,j}(\tA) =  g_{i,j} \TrP \bigodot_{l \neq i,j}a_{l,l}.
\end{equation}
To emphasize, with the above setting when $h=k$, we have the
equalities,
\begin{equation*}\label{eq:GijhkFuncTropicalMultiD}\begin{array}{lll}
   G_{ih,jk}(\tA) = & G_{i,j}(\tA_{h,k}) =
    & g_{i,j}\TrP g_{h,k} \TrP \bigodot_{l \neq i,j,h,k}a_{l,l},
\end{array}
\end{equation*}
and since for this case $g_{h,k} = 0$ then the relation,
$G_{ih,jk}(\tA) = G_{i,j}(\tA)$ is obtained.
\parSpc

\begin{proposition}\label{prop:represOfQuadMatrix}
Let $\tA$ be the symmetric matrix relates to a tropical quadric
$V(\tF)$ in $n$-dimensional space, and let $ a_{i,j} =  \hlf
a_{i,i}\TrP  \hlf  a_{j,j} \TrP \ep_{i,j}$, then
\begin{equation}\label{eq:quadricDecomposition}
    |\tA_{i,j}|= G_{i,j}(\tA) \TrP TropDet(\epMat_{i,j})
\end{equation}
where $\epMat$ is the ``distortion'' matrix of $\tA$.
\end{proposition}
\begin{proof} This assertion is proved by induction over the space
dimension $n$, i.e. the matrix's size minus one.
%\begin{description}
%  \item[Base step for $n=2$:]

\parSpc
\textbf{\emph{Base step for $n=2$:}} Let $\tA$ be the $3 \times 3$
  tropical matrix $(a_{i,j})$, and let
$$
\tA^* = \left[ \begin{array}{lll}
          _{a_{2,2} \TrP a_{3,3} \TrS 2a_{2,3}}&
          _{a_{1,2} \TrP a_{3,3} \TrS a_{1,3} \TrP a_{2,3}}&
          _{a_{1,2} \TrP a_{2,3} \TrS a_{1,3} \TrP a_{2,2}} \\
          _{a_{1,2} \TrP a_{3,3} \TrS a_{1,3} \TrP a_{2,3}} &
          _{a_{1,1} \TrP a_{3,3} \TrS 2a_{1,3}} &
          _{a_{1,1} \TrP a_{2,3} \TrS a_{1,2} \TrP a_{1,3}} \\
          _{a_{1,2} \TrP a_{2,3} \TrS a_{1,3} \TrP a_{2,2}} &
          _{a_{1,1} \TrP a_{2,3} \TrS a_{1,2} \TrP a_{1,3}} &
          _{a_{1,1} \TrP a_{2,2} \TrS 2a_{1,2}}  \\
\end{array} \right].
$$
be its adjoint matrix. Take for instance,
$$ a^*_{3,2} =  \overbrace{\OP a_{1,1} \TrP \hlf a_{2,2} \TrP \hlf a_{3,3}\CP}^{G_{2,3}(\tA)}
\TrP\overbrace{ \OP \ep_{1,1} \TrP \ep_{2,3} \TrS \ep_{1,2} \TrP
\ep_{2,3} \CP}^{|\epMat_{2,3}|}$$
then one can easily verify that it satisfies the required. The
same can be checked similarly for any pair of indices $i,j =
1,2,3$.

%  \item[Induction step $(n-1) \rightarrow n$:]
\textbf{\emph{Induction step $(n-1) \rightarrow n$:}}
  Assuming the correctness
  for $n-1$ we will prove the assertion for $n$.
  Using the minors' decomposition (\ref{diag:tropiclaConvexity}) we can write
\begin{equation*}
  \da_{i,j} = |\tA_{i,j}| =  \bigoplus_{k\neq j} a_{h,k} \TrP
  |\tA_{ih,jk}|, \qquad  h \neq i.
\end{equation*}
Combine this decomposition with the induction assumption and the
description of $a_{h,k}$ in ``distortions'' terms we get,
\begin{equation*}
  |\tA_{i,j}| =  \bigoplus_{k\neq j}
  \OP \hlf a_{h,h}\TrP  \hlf  a_{k,k} \TrP
  \ep_{h,k} \CP \TrP G_{ih,jk}(\tA) \TrP
  |\epMat_{ih,jk}|, \qquad  h \neq i.
\end{equation*}
By the definition of $G_{i,j}(\tA)$, Eq.
\Ref{eq:GijFuncTropicalMultiD}, for any $h \neq i$, $k \neq j$ we
have
$$\OP \hlf a_{h,h}\TrP \hlf a_{k,k}
\CP \TrP G_{ih,jk}(\tA) = G_{i,j}(\tA)$$ and hence,
\begin{equation*}
  |\tA_{i,j}| =  G_{i,j}(\tA) \TrP  \underbrace{\bigoplus_{k\neq j}
  \ep_{h,k} \TrP
  |\epMat_{ih,jk}|}_{= |\epMat_{i,j}|}, \qquad  h \neq i.
\end{equation*}
Clearly the left component is $|\epMat_{i,j}|$.
%\end{description}
\end{proof}

\parSpc

The next two assertions, which refer to some properties of the
``distortion'' matrix $\epMat$, pave the way to the forthcoming
generalization of proposition \ref{prop:conicFullSubDiv} to
quadrics. For this development the following notation is required.

\parSpc
\textbf{\emph{Notation:}} Let $S_n$ be the collection of all the
permutations over the set of indices $N =\{1,\dots,n
  \}$, then $\resPrS{i}{j}$ denotes the sub collection of all permeations $\pi \in
  S_n$ with the fixed correspondence $\pi(i) = j$, namely
\begin{equation}\label{eq:reducePermuatations}
  \resPrS{i}{j} = \{ \pi | \; \pi \in S_n \; s.t. \; \pi(i) = j\}.
\end{equation}

\begin{claim}\label{claim:diffNegativeMatrix}
Let $\epMat = (\ep_{i,j})$ be a $n \times n$ symmetric matrix
where $\ep_{i,i} = 0$  for any $i = 1,\dots, n$ and $\ep_{i,j} <
0$ for any $i \neq j$ then $|\epMat_{i,i}| = 0$ and
$|\epMat_{i,j}| < 0$.
%\begin{enumerate}
%  \item $|\epMat_{i,i}| = 0$,
%  \item $|\epMat_{i,j}| < 0$.
%\end{enumerate}
\end{claim}
\begin{proof} Any minor, $\epMat_{i,i}$, of $\epMat$ is also
a symmetric matrix with zero diagonal. The determinant of any
matrix can be fundamentally stated as $|\epMat| = \bigoplus_{\pi
\in S_n}\bigodot_{i}\ep_{i,\pi(i)}$.
\begin{enumerate}
  \item Since all the matrix's entries are negative except those that  lay
  on the main diagonal which are equal to zero, the single permutation that
  achieves the maximum is the identity (i.e. $\pi(l) =l$), hence
  $|\epMat_{i,i}| = \bigodot_{l \neq i} \ep_{l,\pi(l}) =
    \bigodot_{l \neq i} \ep_{l,l} = 0$.
  \item According to the same considerations as above, without loss of
  generality assume $i<j$, then both the ``new'' $j-1$'th row
  and the ``new'' $i$'th column (of $\epMat_{i,j}$) contain only
  negative elements. As a result for any
  permutation $\pi \in S_n$ we have $\bigodot_{i}\ep_{i,\pi(i)} < 0$ and
  hence use the tropical amount,
  $|\epMat_{i,j}| < 0$.
\end{enumerate}
\end{proof}

\begin{claim}\label{claim:diffpositiveMatrix}
Let $\epMat = (\ep_{i,j})$ be a $n \times n$ symmetric matrix
where $\ep_{i,i} = 0$ for any $i = 1,\dots, n$ and $\ep_{i,j} > 0$
for any $i \neq j$ then
 $|\epMat_{i,i}| \TrP |\epMat_{j,j}| = 2|\epMat_{i,j}|.$
\end{claim}
\begin{proof}
%This assertion is tropically rephrased as
% $$|\epMat_{i,i}| \TrP |\epMat_{j,j}| = |\epMat_{i,j}|\TrP  |\epMat_{i,j}|.$$
%or equivalently
%$
%$ $$\hlf |\epMat_{i,i}| \TrP \hlf |\epMat_{j,j}| =
%|\epMat_{i,j}|.$$
%
Contrarily assume:
 $$2|\epMat_{i,j}| = |\epMat_{i,j}| \TrP |\epMat_{j,i}|
 < |\epMat_{i,i}| \TrP |\epMat_{j,j}|. \qquad (*)$$
Let $N$ be the set of $n$ indices $N =\{1,\dots,n \}$, then any
minor can be written as
  $$|\epMat_{i,j}| = \bigoplus_{\pi \in \resPrS{i}{j}} \bigodot_{k \neq i}
  \ep_{k,\pi(k)}, \qquad (\pi(k) \neq j),$$
where $\resPrS{i}{j}$ is defined according to
\Ref{eq:reducePermuatations}. Recall that, by the assumption
$|\epMat_{i,j}|
> 0$ for any $i,j$.

\parSpc
Let $\pi_{i,i} \in \resPrS{i}{i}$, $\pi_{j,j} \in \resPrS{j}{j}$,
$\pi_{j,i} \in \resPrS{j}{i}$ and $\pi_{i,j} \in \resPrS{i}{j}$ be
the permutations which achieve the maximum for $|\epMat_{i,i}|$,
$|\epMat_{j,j}|$, $|\epMat_{j,i}|$ and $|\epMat_{i,j}|$
correspondingly. Using the above notation, $(*)$ can be rewritten
equivalently as
 $$\overbrace{\underbrace{\bigodot_{k \neq i}
  \ep_{k,\pi_{i,j}(k)}}_{\pi_{i,j}(k) \neq j}}^{fixed \; row}
\TrP
  \overbrace{\underbrace{\bigodot_{k \neq j}
  \ep_{\pi_{j,i}(k),k}}_{\pi_{j,i}(k) \neq i}}^{fixed \; column}
  <
  \overbrace{\underbrace{\bigodot_{k \neq i}
  \ep_{k,\pi_{i,i}(k)}}_{\pi_{i,i}(k) \neq i}}^{I} \TrP
  \overbrace{\underbrace{\bigodot_{k \neq j}
  \ep_{k,\pi_{j,j}(k)}}_{\pi_{j,j}(k) \neq j}}^{II}. $$
Since the permutations in both sizes of the equation are defined
over the same set of indices and are different only in their
restrictions over the indices $i, j$ then the distinction among
their maximum evaluations is derived only from the elements that
correspond to these indices.

\parSpc
Let $k'$ be the index for which $\pi_{i,i}(k') = j$ and let $k''$
be the index for which $\pi_{j,j}(k'') = i$. Set
$\tilde{\pi}_{i,i}$ and $\tilde{\pi}_{j,j}$ to be the results of
the switching $\pi_{i,i}(k') := i$, $\pi_{i,i}(i) := j$ and
$\pi_{j,j}(k'') := j$ $\pi_{j,j}(j) := i$ in the permutations
$\pi_{i,i}$ and $\pi_{j,j}$ respectively. Thus, $\tilde{\pi}_{i,i}
\in \resPrS{i}{j}$ and  $\tilde{\pi}_{j,j} \in \resPrS{j}{i}$.
Since the components are all positive and symmetric (i.e.
$|\epMat_{i,j}| = |\epMat_{j,i}|$), this yields that either the
component that corresponds to the permutation $\tilde{\pi}_{i,i}
\in \resPrS{i}{j}$ is bigger than the one that corresponds to the
preliminary chosen $\pi_{i,j} \in \resPrS{i}{j}$, $\pi_{j,i} \in
\resPrS{j}{i}$ or the second $\tilde{\pi}_{j,j} \in \resPrS{j}{j}$
has this property. Namely, we specify a permutation if
$\resPrS{i}{j}$ which contradicts the preliminary selection of the
permutation which maximized the determinant $|\epMat_{i,j}|$.

\parSpc
Similarly, we can contrarily assume the inequivalence in the
opposite direction
 $|\epMat_{i,j}| \TrP |\epMat_{j,i}|
 > |\epMat_{i,i}| \TrP |\epMat_{j,j}|,$
 to obtain a contradiction as well. This completes the proof of the
 required.
\end{proof}

\parSpc
Using the above preparation we can now generalize proposition
\ref{prop:conicFullSubDiv} to quadrics.
\begin{theorem}\label{prop:conicFullSubDivMulti} Let $\tF$ be a
 quadratic tropical polynomial and let $\tS_\tF$ be the corresponding induced
 subdivision of its Newton polytope $\Polygon_\tF$.
 The dual subdivision $S_{\dtF}$ (induced by $\tF^*$) of
\begin{enumerate}
  \item a subdivision $\tS_\tF$ that is \full, is \full,
  \item a subdivision $\tS_\tF$ that is \emty, is \emty.
\end{enumerate}
\end{theorem}
\begin{proof}
Let $\tA$ be the tropical symmetric matrix that characterizes
$\tF$ and let $\dtA$ be the corresponding matrix that describes
$\dtF$ -- the dual of $\tF$. The subdivision $\tS_{\tF}$ is
induced directly from the translation of $\tA$ to $\Polygon_\tF$
\Ref{diag:genMatrixToPolytop} and similarly from the dual (i.e.
$\dtA$ to ${\Polygon}_{\dtF}$). To prove these two assertions we
will assisted  by Proposition \ref{prop:represOfQuadMatrix} which
is  based on the entries' representation
(\ref{eq:conicConvexityInTropicalMultiD}) in ``distortions''
terms.

\parSpc
First, in both cases, due to the property of $G(\tA)$ (derived
from Def. \ref{eq:GijFuncTropicalMultiD}) which determines
$$ G_{i,i}(\tA) + G_{j,j}(\tA)  = 2G_{i,j}(\tA)$$
when observing the convexity of any triple,
$$
\begin{array}{ccccc}
  |\tA_{i,i}| & \longleftrightarrow & |\tA_{i,j}| & \longleftrightarrow & |\tA_{j,j}|, \\
%
%  || &  & || & & || \\
%  G_{i,i}(\tA) |\epMat_{i,i}| & \longleftrightarrow & |\tA_{i,j}| & \longleftrightarrow & |\tA_{j,j}| \\
\end{array}
$$
the tropical multipliers, $G_{i,j}(\tA)$, can be neglected. Thus,
convexity examination can be based only on the ``distortion''
values $\ep_{i,j}$.
\begin{enumerate}
%%%--------------------------
  \item
  By the setting (\ref{eq:conicConvexityInTropicalMultiD}), a subdivision
  $\tS_\tF$ is \full \; if and only if
  $ \ep_{i,j} < 0 $ for all $i \neq j$.
  As a result of Claim \ref{claim:diffNegativeMatrix} the
  following translation is received:

$$ a^*_{i,j} = |\tA_{i,j}| \longmapsto |\epMat_{i,j}|   \left\{
\begin{array}{cc}
    = 0 & i = j, \\ [2mm]
    < 0 & i \neq j.
\end{array}
\right.
$$
%$$ \begin{array}{ccccc}
%    a^*_{i,i} & \longleftrightarrow & a^*_{i,j} & \longleftrightarrow & a^*_{j,j} \\
% || & & || & & || \\
% |\tA_{i,i}| & \longleftrightarrow & |\tA_{i,j}|& \longleftrightarrow & |\tA_{j,j}| \\
%  \downarrow & & \downarrow & & \downarrow \\
%  |\epMat_{i,i}| & \longleftrightarrow &|\epMat_{i,j}| & \longleftrightarrow &
%  |\epMat_{j,j}| \\
%   || & & \wedge & & || \\
% 0 & \longleftrightarrow & 0 & \longleftrightarrow & 0 \\
%\end{array}$$

As one can easily verify,  any node which is not vertex (i.e. is
corresponding to indices $i \neq j$) has a negative value and
hence, it appears in the lower part $\lbCH(\Polygon_{\dNuF})$.
This means the it is expressed in the compatible subdivision,
$\tS_{\dtF}$, of the dual polytope $\Polygon_{\dtF}$. This
property is valid for any $i \neq j$, or equivalently for any node
which is not vertex, hence, the subdivision is \full.
%%%--------------------------
  \item
 A subdivision $\tS_\tF$ is \emty \;  if and only if
  $ \ep _{i,j} \geq 0 $ for $i,j=1,\dots, n+1$ while $\ep_{i,i} = 0$.
  The same convexity consideration as has been used in the previous
  part is also standing in this case. Using Claim
  \ref{claim:diffpositiveMatrix} we have the equality
  $2|\epMat_{i,j}| = |\epMat_{i,i}| \TrP |\epMat_{j,j}|$. Thus,
  despite $\da_{i,j}$ being  placed on the ``lower part'' of the convex
  hull, $\lbCH(\Polygon_{\dNuF})$, it is not one of the polytopes  vertices. This  means
  that any middle node does not appear in the subdivision $S_{\dtF}$
  which is
  induced by $\lbCH(\Polygon_{\dNuF})$. The net result is that
  only the
  vertices of Newton Polytope appear in the subdivision, hence the subdivision is \emty.
\end{enumerate}
\end{proof}

%~~~~~~~~~~~~~~~~~~~~~~~~~~~ SUB-SECT ~~~~~~~~~~~~~~~~~~~~~~~~~~~~~~~~~~~
\secSpc
\section{Regularity in Duality Sense}\label{sec:RegularityInDualitySense}
Up to this point, all our development was done for families of
objects with respect to their dual, where the referred issues were
either the resulting geometric objects or equivalently the
subdivisions determining these objects. Here, we intend to advance
these idea and to specify additional properties of the
subdivision's generators, namely the lifting functions. Recall
that lifting functions are defined by the coefficients of the
tropical polynomials. For this purpose we will specify lifting
functions $\nuF$ which have a ``cooperative'' behavior with
respect to applying duality twice. The corresponding tropical
curve, $V(\tF)$, of such functions will be denoted as
\emph{Regular Tropical Curves}. The formal meaning is described
next.
\parSpc
\textbf{\emph{Notation:}} Let $\tF$ be a tropical polynomial
defined by the matrix $\tA$ and $\dtF$ its dual (determined by
$\dtA$), then $\ddtF$ denotes the dual of $\dtF$ (i.e. the result
of applying twice the duality to $\tF$) and $\ddtA$ is its
corresponding matrix.

\parSpc
\begin{definition} %{\bf (Regular of Tropical Duality)}
Let $V(\tF)$ be a tropical curve where $\tF$  is described by the
symmetric matrix $\tA$. The curve will be called \textbf{Regular}
in duality sense if
$$a^{**}_{i,j} = \lC{\tA} \TrP a_{i,j}, \qquad \forall i,j = 1, \dots, n,$$
where $a^{**}_{i,j}$ is an entry of $(\dtA)^{*}$ and $\lC{\tA}$ is
a fixed value which is called the lifting constant.
\end{definition}

\parSpc
By the above, the corresponding convex hulls of the lifting
functions $\nuF$ and $\ddNuF$ are same up to lift by constant.
Thus, it is clear that the regularity in this sense means that for
the two tropical polynomials $\tF$ and $\ddtF$ the induced
subdivisions $S_{\tF}$ and $S_{\ddtF}$ are identical. In this case
we may also say that the polynomial $\tF$ and the lifting function
$\nuF$ are regular. The two following claims rephrase the
regularity in term of ``distortion'' values, these are a
preparation for the main statement.
\begin{claim}\label{claim:regualrityOfnegtiveByNorm}
 Let $\tF$ be a
 quadratic tropical polynomial and $\ddtF$ the result of
 applying twice the duality to $\tF$, and let $\tA$, $\ddtA$ be their
 characteristic  matrices having the ``distortion'' matrices
  $\epMat$ and $\ddEpMat$
 correspondingly. In case $\ep_{i,j} < 0$ for all $i \neq j$,
  $V(\tF)$ is regular if and only if  $\epMat =
 \epMat^{**}$.
\end{claim}

\begin{proof} This assertion is directly derived for Claim
\ref{claim:diffNegativeMatrix}. For this case $\ep_{i,i}=
\dEp_{i,i} = \ddEp_{i,i} = 0$. Thus, the lifting constant is
$\lC{\tA}=0$ and hence $a_{i,i}= \dda_{i,i}$. In addition, since
$\lC{\tA}=0$ we have $\ep_{i,j}= \ddEp_{i,j}$ for any $i,j
=1,\dots,n$ then, $\dda_{i,j} = \hlf \dda_{i,i} \TrP \hlf
\dda_{j,j} \TrP \ddEp_{i,j} = \hlf a_{i,i} \TrP \hlf a_{j,j} \TrP
\ep_{i,j} = a_{i,j}$.
\end{proof}

\begin{claim}\label{claim:MinorsOfnegtiveByNorm}
Let $\epMat = (\ep_{i,j})$ be a $n \times n$ symmetric matrix
where $\ep_{i,i} = 0$ for any $i = 1,\dots, n$ and $\ep_{i,j} < 0$
for any $i \neq j$ then
$$ |\epMat_{i,j}|=\ep_{j,i} \TrS
\OP \bigoplus_{k \neq i,j}\ep_{j,k}\TrP \bigoplus_{h \neq
i,j}\ep_{h,i}\CP   \; = \;  \left\{ \begin{array}{cc}
  |\epMat_{i,i}| = 0, & i= j, \\ [2mm]
  |\epMat_{i,j}| < 0, & i \neq j.
\end{array}\right.$$
%$$ |\epMat_{i,j}|=\ep_{i,j} \TrS
%\left(\bigoplus_{k \neq i,j}\ep_{j,k}\TrP \bigoplus_{h \neq
%i,j}\ep_{h,i}\right).$$
%
\end{claim}
\noindent \textbf{\emph{Remark:}} This assertion is also valid
without the symmetry restriction.

\begin{proof}
A careful combinatoric observation of the matrix decomposition
yields the result. Without loss of generality assume $j>i$, then
the minor $\epMat_{i,j}$ is of the form:
$$\epMat_{i,j}
=\begin{array}{c|cccc}
      &   & i  &   &   \\ \hline
      & 0 & |  &   &   \\
      &   & |  &   &   \\
  j-1 & - & \ep_{j,i} & - & -  \\
      &   & |  &   & 0,
\end{array}$$
where, besides the $i$'th column and $(j-1)$'th row, any row or
column contains a zero entry. The tropical determinant ``selects''
the maximal product that is achieved by some permutation which
simultaneously choses one element from each row and each column.
Naturally, in this case since all the non-diagonal entries are $<
0$, the aim is to chose the maximal number of zero entries.

\parSpc However, such choosing should be completed with
entries of the $i$'th column and $(j-1)$'th row. One option is to
select their intersection (i.e. the entry $\ep_{j,i}$), while the
other option is a separated selection, one for the column and one
for the row, and then taking their tropical product.  Eventually,
we take the ``best'' result that can be achieved from both
options. Note that in these two scenarios the zero entries don't
play a role. Formulating the above we get,
$$ |\epMat_{i,j}|=
\bigoplus_{k \neq j}\ep_{j,k}\TrP|\epMat_{ij,jk}| = \ep_{j,i} \TrS
\OP  \bigoplus_{k \neq i,j}\ep_{j,k}\TrP \bigoplus_{h \neq i,j}
\ep_{h,i} \CP  ,
$$
and the required is attained. Moreover, since this is a
combination of values $\ep_{i,j} < 0 $ then for $i \neq j$
$|\tE_{i,j}| < 0$. In case $i=j$, each column and each row of
$|\tE_{i,i}|$ contains a zero entry, which is the maximal chosen,
thus $|\tE_{i,i}| = 0$.
\end{proof}
The next two theorems give the explicit conditions for determining
whether or not a curve of the previous families is regular.

\begin{theorem}\label{prop:conicFullSubDivMulti} Let $\tF$ be a
 quadratic tropical polynomial, $\tA$ its characteristic matrix and let
 $\epMat$ be the  ``distortion'' matrix of $\tA$.
 Assume $\ep_{i,j}<0$ for all $i,j=1,\dots,n$ and $\ep_{i,i} = 0$,
 a sufficient and necessary condition for $V(\tF)$ to be regular
 is
 $$
 \ep_{j,i} >
\bigoplus_{k \neq i,j}\ep_{j,k}\TrP \bigoplus_{h \neq i,j}
\ep_{h,i} \qquad (*)
 $$
  for any $i,j = 1,\dots,n$.
\end{theorem}
\begin{proof}
Using Claim \ref{claim:MinorsOfnegtiveByNorm}, due to the
assumptions, $(*)$ can be  stated equally as $\dEp_{i,j}=
\ep_{i,j}$. The sufficiency part is a direct result of Claim
\ref{claim:MinorsOfnegtiveByNorm}, for this case we have
$\dEp_{i,j}= |\epMat_{i,j}| =  \ep_{i,j}$ and thus $\ddEp_{i,j}=
|(\dEpMat)_{i,j}| =  \dEp_{i,j}= \ep_{i,j}$ which yields $\epMat =
\ddEpMat$. Then apply Claim \ref{claim:regualrityOfnegtiveByNorm}
to obtain the required.
\parSpc
For the necessary part, contrarily assume there exist indices
$i,j$ such that $(*)$ is not satisfied, namely $\dEp_{i,j} >
\ep_{i,j}$, and the equivalence $\epMat = \ddEpMat$ is still
valid. Apply Claim \ref{claim:MinorsOfnegtiveByNorm} for $\dEpMat$
to obtain:
$$ \ddEp_{i,j} = |\ddEpMat_{i,j}|=\dEp_{i,j} \TrS
\left(\bigoplus_{k \neq i,j}\dEp_{j,k}\TrP \bigoplus_{h \neq
i,j}\dEp_{h,i}\right) > \ep_{i,j}.$$ That is a contradiction to
Claim \ref{claim:regualrityOfnegtiveByNorm} which determines the
equivalence $ \ep_{i,j}=\ddEp_{i,j}$.
\end{proof}

\begin{theorem}\label{prop:conicFullSubDivMulti} Let $\tF$ be a
 quadratic tropical polynomial, $\tA$ its $n \times n$
 characteristic  matrix and let $\epMat$ be the ``distortion'' matrix of
 $\tA$. In case $\ep_{i,j}>0$ for all $i,j=1,\dots,n$, $V(\tF)$ is never regular.
\end{theorem}
\begin{proof}
Using claim \ref{claim:diffpositiveMatrix} we obtain that
$\hlf |\epMat_{i,i}| \TrP \hlf |\epMat_{j,j}| = |\epMat_{i,j}|,$
for any $i,j$. As a result,  the ``distortion'' matrix of
$\dEpMat$, is a complete zero matrix $\tZ$. Since all of its
entries are equal, we get $|\tZ_{i,j}| = -\infty$. According to
the definition of the ``distortion'' values the equality
$|\epMat_{i,j}| = \lC{\epMat} \TrP |\trop{Z}_{i,j}| = \lC{\tA}
\TrP -\infty = -\infty$ is satisfied.  Using the same
consideration again we obtain $|\tA_{i,j}|= -\infty $, and hence
$\da_{i,j} = -\infty$ for all $i,j$. Thus any combination of these
entries has the same evaluation and in particular $|(\dtA)_{i,j}|
= \dda_{i,j} = -\infty \dir a_{i,j} \neq \lC{\tA} \TrP \dda_{i,j}
$. Hence, $V(\tF)$ is non-regular.
\end{proof}
%******************************* SECT 5 *********************************
%\secSpc
%\section{Conclusions}\label{sec:Conclusions}
%
%******************************* SECT  *********************************
\section*{\center Appendix}
\begin{appendix}
\secSpc
\section{Duality of Algebraic Varieties}
In this appendix we will overview the classical algebraic duality
in polynomial equations \cite{griffiths78principles,hartshorne79},
this review will include the explicit construction (i.e. an
algorithm) of the polynomial equation which describes the dual
curves. We open by the formal definition of duality that occurred
in the projective space, which is used later to define the duality
over affine spaces.

%~~~~~~~~~~~~~~~~~~~~~~~~~~~ SUB-SECT ~~~~~~~~~~~~~~~~~~~~~~~~~~~~~~~~~~~
\subSecSpc
\subsection{Projective Dual Curve}

%Let $\curve \subset \Proj^2$, the set of lines in $\Proj^2$ are
%regraded as another projective space $(\Proj^*)^2$, by taking
%$(a_0:a_1:a_2)$ as homogeneous coordinates of the line $L: a_0 +
%x_0 a_1 x_1 + a_2 x_2 =0$
%The \emph{Dual Curve}
% ----------------------------
Let  $\Proj^2_\Fld = Proj \; \Fld[T_0,T_1,T_2]$ be the projective
plane, a curve $L$ on $\Proj^2_\Fld$ is said to be a \emph{line},
if $Deg(L) = 1$ i.e. $L = V(a_0 T_0 + a_1 T_1 + a_2 T_2)$ for some
$(a_0,a_1,a_2) \in (\Fld^*)^3$. Clearly
$$
V(a_0 T_0 + a_1 T_1 + a_2 T_2) = V(b_0 T_0 + b_1 T_1 + b_2 T_2)
$$
$$
  \Longleftrightarrow  \qquad a_i = \lambda b_i, \qquad \; \; i = 0,1,2, \;
  \;  \lambda \in \Fld^*.
$$
Therefore a bijection
$$
\phi: \{lines \; on \; \Proj^2_\Fld \} \; \; \longrightarrow \; \;
\Proj^2_\Fld (\Fld)
$$
is obtained and defined by
$$\phi(V(a_0 T_0 + a_1 T_1 + a_2 T_2)) = (a_0:a_1:a_2).$$
Namely, $(a_0:a_1:a_2)$ are taken as homogeneous coordinates of
the line $L: a_0 T_0 + a_1 T_1 + a_2 T_2 =0$.

\parSpc
Let $\curve =V(F)$ be a curve with $Deg(\curve) = Deg(F) =r$,
where $F$ is homogeneous polynomial of degree $r \geq 2$. For a
nonsingular point $p$ of $\curve$, let $l_p$ denote the tangent to
$\curve$ at $p$, equivalently the tangent space $T_p$ at the point
$p$. In case $p = (p_0, p_1, p_2 )$ and $(\partial F / \partial
T_i) (p) = \partial_i F(p_0, p_1, p_2 )$ is written as
$(\partial_i F)(p)$ for $i=0,1,2$, and
$$l_p = V\left(\sum_{i=0}^2(\partial_i F)(p)T_i\right),$$
i.e. $ \phi (l_p) = ((\partial_0 F)(p):(\partial_1
F)(p):(\partial_2 F)(p))$. From the rational map $\theta_\curve$
defined by % $(\partial_i F)(p)$ is denoted by $\theta_\curve$
$\theta_\curve(p)=((\partial_0 F)(p):(\partial_1 F)(p):(\partial_2
F)(p))$ the relation  $\theta_\curve(p) = \phi(l_p)$ between the
maps' images is trivially obtained.

\begin{definition}
 $\curve^*= \theta_\curve(\curve)$ is said to be the dual curve of
 $\curve$. In case $\curve^*$ is expressed as $V(F^*)$ by an
 irreducible homogeneous polynomial  $F^*$, then  $F^*$ is said to
 be the dual homogeneous polynomial of $F$.
\end{definition}

%~~~~~~~~~~~~~~~~~~~~~~~~~~~ SUB-SECT ~~~~~~~~~~~~~~~~~~~~~~~~~~~~~~~~~~~
\subSecSpc
\subsubsection{Construction of Dual Curve in $\Proj^2$}

In the pervious section we outlined the formal definition of the
dual curves with the following relation
\begin{equation*}
  \begin{CD}
    \curve                        @>\theta_\curve>>             \curve^*   \\
      @|                                   @                  |  \\
      V(F)                                 @> >>
      V(F^*).
 \end{CD}
\end{equation*}
The next goal is, to attain the description of a dual curve
$\curve^*$ according to the primal curve $\curve$. Namely, to
write down the explicit equation of $F^*$ where its coefficients
are described in terms of $c_{\om}$'s (the coefficients of $F$).
For our concern we may assume that $F$ is a Laurent polynomial. A
number of preparatory steps are needed in order to pave the way
for the easier construction. In order to simplify these geometric
matters one can refer to the projective plane as Euclidian
3-dimensional space (in this case the coordinates are signed by
$(t_0,t_1,t_2)$) with all the required restrictions.

\parSpc
%\begin{description}
%\item[\emph{Step 1:}]
\textbf{\emph{Step 1:}}
% -------- step 1 ----------
In order to define the tangent space $T_p(\curve)$ (line in this
case), the \emph{Gradient} $\nabla F$ is found and denoted by :
$$ \nabla F (p) = ((\partial_0 F)(p),(\partial_1
F)(p),(\partial_2 F)(p)) = (\al_0, \al_1, \al_2 ).
$$
The derivatives provide the direction of the normal to the tangent
space at a fixed point $p^0=({p}_0^0, {p}_1^0, {p}_2^0)$ that
satisfies $F$. Taking $p^0$ as varied along, the curve of the
complete family of the tangent spaces  $\tanBun_\curve$ (i.e.
tangent bundle) is obtained.

\parSpc
Recall that homogeneous polynomial of degree $n$ satisfying  the
relation
\begin{equation}\label{eq:homoPoly}
t_0(\partial_0 F) + t_1(\partial_1 F) + t_2(\partial_2 F)
 = n F.
\end{equation}
Hence, when $ F = 0$ the equations of the tangent lines can be
written as
\begin{equation}\label{eq:tangentLine}
  \al_0 t_0 + \al_1 t_1 + \al_2 t_2 = 0.
\end{equation}
\parSpc
%\item[\emph{Step 2:}]
\textbf{\emph{Step 2:}}
% -------- step 2 ----------
A parameterization of the tangent bundle $\tanBun_\curve$ along
$\curve$ can be obtained due to the substitution
\begin{equation}
 \label{eq:subTangentBundle}
  t_0  \mapsto -(\frac{\al_1 u + \al_2 v}{\al_0}),
  \qquad t_1 \mapsto u, \qquad t_2 \mapsto v,
\end{equation}
in $F$, and by the homogeneity of $F$ we have
\begin{equation}\label{eq:glines}
  \tilde{F} = F(-(\al_1 u + \al_2 v), \al_0 u, \al_0 v) = 0.
\end{equation}
Note that this polynomial equation is homogeneous as well and of
the same degree in its \emph{five} variables (i.e. $\tilde{F} \in
\Fld[\al_0,\al_1,\al_2,u,v]_n$). However, it can also be regarded
as polynomial only in $u$ and $v$ (i.e. $\tilde{F} \in
\Fld[\al_0,\al_1,\al_2][u,v]_n$) where its coefficients are
determined by the $\al_i$'s
\begin{equation}\label{eq:Fne}
  \tilde{F}_{(\al_0, \al_1, \al_2)}(u, v) = \sum_{i =0}^n a'_{i}
  (\al_0, \al_1, \al_2){u}^i {v}^{n-i} = 0 \; .
\end{equation}
%
%\item[\emph{Step 3:}]
\textbf{\emph{Step 3:}}
% -------- step 3 ----------
At last, the required dual curve $\curve^*$ is the projection of
$\tanBun_\curve$ into the space which spanned  by $\al_i$'s.
Algebraically, this is equivalent to the elimination of $u$ and
$v$ from $\tilde{F}$. Here the homogeneity of $\tilde{F}$ and the
property of being a Laurent polynomial becomes useful .
\begin{lemma}
Let $F$ be a non-degenerate homogeneous polynomial in $t_i$ then
$$F(p) = 0 \qquad \Longleftrightarrow \qquad \frac{\partial
F}{\partial t_i}(p) = 0, \forall i  . $$
\end{lemma}
\noindent The lemma is directly derived from the property
(\ref{eq:homoPoly}) and the assumption that $F$ is not a
degenerate polynomial.

\parSpc
Suppose that $A = \{m_1,\dots,m_s\} \subset \Int^n$ and $A$
generates $\Int^n$. Let $L(A)$ be the set of Laurent polynomials
with exponent vectors in $A$, i.e.,
$$ L(A) =\{ a_1 t^{m_1} + \dots + a_s t^{m_s}: a_i \in \Fld \}$$
where $t^m = t^a_1\dots t^a_n$ for $m = (a_1,\dots,a_n ) \in
\Int^n$. Given $n+1$ Laurent polynomials $f_0,\dots,f_n \in L(A)$,
their A-resultant
$$
Res_A (f_0,\dots,f_n)
$$
is a polynomial in the coefficients of the $f_i$.
\begin{proposition}
The vanishing of $Res_A (f_0,\dots,f_n)$ is necessary and
sufficient condition for the equations $f_0 = \dots = f_n = 0$ to
have a solution (see [\cite{Gelfand94}, Prop. 2.1]).
\end{proposition}

However, one must be careful where the solution lies. The $f_i$
are defined initially on the torus $(\Comp^*_m)^n$, but the
definition of the generalized corresponding  projective toric
variety $Y_A$ shows that the equation $f_i = 0$ makes sense on
$Y_A$ . Then one can prove that
$$Res_A (f_0 , \dots , f_n ) = 0 \; \Longleftrightarrow \; f_0 = \dots = f_n = 0$$
have a solution in $Y_A$. This stays true for any $A \subset
\Int^n$ finite or not, especially while $A$ is bounded from below
which is enough for our purpose.

%\parSpc
%The following is  a direct result of the above:
%
\begin{proposition}
Let $\tilde{F}$ be a homogeneous Laurent polynomial then
$$\tilde{F} = 0 \qquad \Longleftrightarrow \qquad Res ( \frac{\partial
\tilde{F}}{\partial u},\frac{\partial \tilde{F}}{\partial v} ) =
0.$$
\end{proposition}
\noindent As a result the dual curve $\curve^*$ is described by
$V(Res(\al_0,\al_1,\al_2))$.
%\end{description}

%~~~~~~~~~~~~~~~~~~~~~~~~~~~ SUB-SECT ~~~~~~~~~~~~~~~~~~~~~~~~~~~~~~~~~~~
\subSecSpc
\subsubsection{Affine Dual Curve}

The dual curve had been defined for the case of projective spaces,
next we intend to show the compatible duality over affine spaces.
In general, two steps are required, the first to embed the affine
curve defined by $f$ in $\Proj^2$, and the second to ``extract''
the dual from $\Proj^2$ after it has been determined. This can be
outlined via the following diagram
\begin{equation*}
  \begin{CD}
    V(F)                        @>Dual >>             V(F^*)   \\
      @A AA                                   @                      V extract VV \\
      V(f)                                 @> >> V(f^*).
 \end{CD}
\end{equation*}
Note that in this case the dual of a line is defined as follows:
\begin{equation}\label{eq:affDual}
  l: t_2 = a_1 t_1 + a_0 \; \longmapsto \;  (a_1,-a_0),
\end{equation}
and the tangent in a given point $p^0= (p^0_0,p^0_1)$ is
\begin{equation}\label{eq:affTang}
  l: t_2 = \left(\frac{\partial f }{ \partial t_1} \right)(p^0) t_1 + p^0_1 -
  \left(\frac{\partial f }{ \partial t_1}\right)(p^0)
  p^0_0.
\end{equation}

\parSpc
%\begin{description}
\textbf{\emph{Step 1:}}
%\item[\emph{Step 1:}]
% -------- step 1 ----------
Naively embed ($(t_1,t_2) \mapsto  (1,t_1,t_2)$) the curve
$\curve$ defined by $f$ in the projective plane to obtain the
curve described by the homogeneous polynomial $F$ (which is
irreducible if $f$ dose) where $F|_{t_0 = 1 } = f$. Then find
$F^*.$

\parSpc
%\item[\emph{Step 2:}]
\textbf{\emph{Step 2:}}
% -------- step 2 ----------
The ``extraction'' of the dual $f^*$ is based on the tangent space
as constructed in (\ref{eq:tangentLine}) and the restriction of
$t_0$ to 1. In this case one has $\al_0 = -(\al_1 t_1 + \al_2
t_2)$. Rewriting the dual transform (\ref{eq:affDual}) in terms of
$F$ and using the \emph{Implicit Function Theorem}, the image of a
tangent to the curve defined by $f$ is
\begin{equation}\label{eq:affDualTrans}
 a_1 = -\frac{\al_1}{\al_2} \; , \qquad  a_0 % =t_2-\frac{\al_1}{\al_0}t_1
  = \frac{\al_2 t_2 + \al_1 t_1}{\al_2}.
\end{equation}
Composing this with the restriction of $t_0 = 1$ we obtain
\begin{equation}\label{eq:F2final}
 a_1 = -\frac{\al_1}{\al_2} \; , \qquad  a_0 % =t_2-\frac{\al_1}{\al_0}t_1
  = -\frac{\al_0}{\al_2}.
\end{equation}
Let $c =\al_2$ where $0 \neq c \in \Fld$, so that $\al_1 = - c
a_1$ and $\al_0 = -c a_0$.
 The substitution of these values provides the transform of the original curve
$f(t_1, t_2) = 0$ :
$$ Res(\al_0, \al_1, \al_2) = Res ( -c a_0, - c a_1  , c) = $$
$$c ^m Res ( - a_0,  - a_1  , 1) = 0 \; \; \Rightarrow $$
$$Res ( - a_0,  - a_1  , 1) = 0 \; ,$$
where $c \not = 0$ and $m \leq n(n-1)$.
%\end{description}
The net result is, the dual curve for this case, is described by
$V(Res( - a_0, -a_1,1))$.
%~~~~~~~~~~~~~~~~~~~~~~~~~~~ SUB-SECT ~~~~~~~~~~~~~~~~~~~~~~~~~~~~~~~~~~~
\subSecSpc
\subsection{Duality of Quadrics in
Matrices' Notion}\label{subSec:DualityOfMultidimensionalQuadrics}

In this section we apply the complete dual transformation to the
general form of conics, then using this base we will construct the
duality of quadrics in matrices' notion.

\parSpc
Let $f$ be a polynomial of degree 2 defined as
\begin{equation}\label{eq:conicMatForm}
 f( \bar{z}) = \bar{z} C \bar{z}^T =0
\end{equation}
where $ \bar{z} =  [z_1 ,z_2 , 1 ]$  and
\begin{equation}\label{eq:conicMat}
  C = \left[
\begin{array}{ccc}
  c_1 & c_4 & c_5 \\
  c_4 & c_2 & c_6 \\
  c_5 & c_6 & c_3
\end{array} \right],
\end{equation}
the obtained result of applying the duality is $f^*=
\bar{a}C^*\bar{a}^T$ where $\bar{a} = [a_1 ,a_2 , 1 ]$ and

\begin{equation}\label{eq:conicDualMat}
C^*= \left[ \begin{array}{lll} (c_{6}^{2}-c_{3}c_{2}) & (
c_{6}c_{5} -c_{4}c_{3}) & ( c_{5}c_{2}-c_{4}c_{6}) \\
( c_{6}c_{5} - c_{4}c_{3}) & ( c_{5}^{2}-c_{1}c_{3})  & (
c_{4}c_{5}- c_{1}c_{6}) \\ ( c_{5}c_{2}-c_{4}c_{6}) & ( c_{4}c_{5}
- c_{1}c_{6}) & ( c_{4}^{2}-c_{1}c_{2})
\end{array} \right].
\end{equation}
\parSpc
Since for the case of conics there is a preserving of degrees
\cite{walker78algebraic}, the duality induced the mapping of
matrices
$$ C_{3 \times 3} \; \longmapsto \; C_{3 \times 3}^* $$
and vise versa.
Reapplying the dual transform, one can easily check that
\begin{equation}\label{eq:conicDualDual}
(C^*)^* = -Det(C) C.
\end{equation}
Considering the conics which are varieties of a given describing
equation of the form $f = \bar{z}C \bar{z}^T =0$, a corresponding
duality occurs in the space of the compatible  matrices $C_{3
\times 3}$. This relation can generalized to quadrics embedded in
space of any dimension.

\parSpc
Reobserving $C^*$ as appears in equation (\ref{eq:conicDualMat})
one can easily verify that
$ C^* = -Adj(C)$
is the standard adjoint matrix of $C$. Additionally, since only
$V(f^*)$ are considered, any multiplier $d$ of $C^*$ is acceptable
without harming the result, particularly for the setting $d :=
-1/Det(C)$. Composing this specific setting with the pervious
observation yields,
\begin{equation}\label{eq:conicByInverseMat}
C^* = -Adj(C) \; \longleftrightarrow \; -Adj(C)(-\frac{1}{Det(C)})
= C^{-1}.
\end{equation}
Clearly the above coincides with both the direct calculation as
appears in (\ref{eq:conicDualDual}) and the applying of the
duality twice
$$(C^*)^* \; \longleftrightarrow \; (C^{-1})^{-1}=C.$$

\parSpc
Next, we will construct directly the dual of quadrics in terms of
matrices, this can be done due to their special properties. Let
\begin{equation}\label{eq:generalQuadric}
  f(\brZ) = \brZ C \brZ^T=0
\end{equation}
be a polynomial of degree 2 in $\Comp^n$ where $C$ is
$(n+1)\times(n+1)$ symmetric matrix and $\brZ = (z_1, \dots, z_n ,
1)$. Its \emph{Gradient} $\nabla f = 2C\brZ$ defines the cotangent
bundle $\coTanBun_{V(f)}$ of the corresponding manifold $V(f)$.
Hence, the obtained map for each point $\brZ \in V(f)$ is
$$ 2C\brZ^T|_{z_0} \; \longmapsto \; (\brZ_0^*)^T.$$
Assuming $C$ is an inventible matrix and since the above is valid
for any point we have
\begin{equation*}
 \brZ^T = \frac{1}{2}C^{-1}(\brZ^*)^T.
\end{equation*}
Substitute this into equation (\ref{eq:generalQuadric}) to get
\begin{equation*}
  \left[ \frac{1}{2}C^{-1}(\brZ^*)^T \right]^T C \left[\frac{1}{2}C^{-1}(\brZ^*)^T
  \right]=0.
\end{equation*}
Developing the above via,
\begin{equation*}
 \left[ \frac{1}{2}C^{-1}(\brZ^*)^T \right]^T C \left[\frac{1}{2}C^{-1}(\brZ^*)^T \right]=
  \frac{1}{4}\left[(\brZ^*)^T\right]^T \left[C^{-1}\right]^T C
\left[C^{-1}\right]\left[(\brZ^*)^T \right] = 0
\end{equation*}
yields
\begin{equation*}
  \frac{1}{4}\brZ^*\left[C^{-1}\right]^T (\brZ^*)^T =0.
\end{equation*}
Since only the varieties of this equation are concerned we remain
with
\begin{equation}\label{eq:generalQuadricsDual}
  \brZ^*\left[C^{-1}\right]^T (\brZ^*)^T
  =0.
\end{equation}
Using the matrices identity $(C^{-1})^T = (C^{T})^{-1}$ and the
fact that $C$ is symmetric, this result coincides with the basic
demand that the dual of a dual is the source.
\end{appendix}
%******************************* Reference *********************************
%\bibliographystyle{abbrv}
%\bibliography{../../bib/dfz}

\end{document}

%% file: tropMacro.tex
\usepackage{amsfonts,amssymb,amsmath,amscd,amsthm}
\usepackage[mathscr]{eucal}
\usepackage{mathrsfs}
\usepackage{oldgerm,units}
\usepackage{wrapfig,epsfig}

\usepackage{ifthen}
\input diagram.tex

\newtheorem{theorem}{Theorem}[section]
\newtheorem{proposition}[theorem]{Proposition}
\newtheorem{definition}[theorem]{Definition}

\newtheorem{claim}[theorem]{Claim}
\newtheorem{lemma}[theorem]{Lemma}

\newtheorem{corollary}[theorem]{Corollary}

\newtheorem{comment}[theorem]{Comment}
\newtheorem{example}[theorem]{Example}
\newtheorem{remark}[theorem]{Remark}
\newtheorem{conclusion}[theorem]{Conclusion}

\newtheorem{observation}[theorem]{Observation}
\newtheorem{terminology}[theorem]{Tterminology}

\newcommand{\Ref}[1]{(\ref{#1})}
\newcommand{\Comp}{\mathbb C}
\newcommand{\Real}{\mathbb R}
\newcommand{\Rati}{\mathbb Q}
\newcommand{\Int}{\mathbb Z}

\newcommand{\Proj}{\mathbb P}

\newcommand{\Fld}{\mathbb K}

\newcommand{\trop}[1]{\mathcal{#1}}
\newcommand{\Amb}{\trop{A}}
\newcommand{\nAmb}{\trop{N}}
\newcommand{\tA}{\trop{A}}

\newcommand{\tE}{\trop{E}}
\newcommand{\tF}{\trop{F}}

\newcommand{\tS}{\trop{S}}

\newcommand{\tZ}{\trop{Z}}

\newcommand{\To}{\longrightarrow }
\newcommand{\Relate}{\rightsquigarrow }

\newcommand{\dir}{\Longrightarrow}

\newcommand{\FigWidth}{2.5}

\newcommand{\tUniS}{-\infty}

\newcommand{\du}[1]{#1^*}
\newcommand{\ddu}[1]{#1^{**}}
\newcommand{\df}{\du{f}}
\newcommand{\dtA}{\du{\tA}}
\newcommand{\dtF}{\du{\tF}}

\newcommand{\ddtA}{\ddu{\tA}}
\newcommand{\ddtF}{\ddu{\tF}}

\newcommand{\om}{\omega}

\newcommand{\ep}{\epsilon}
\newcommand{\al}{\alpha}

\newcommand{\prS}{S_n}
\newcommand{\resPrS}[2]{\prS|_{#1 \mapsto #2}}

\def \full {maximal in nodes}
\def \emty {minimal in nodes}

\newcommand{\dNu}{\du{\nu}}
\newcommand{\nuf}{\nu_f}
\newcommand{\nuDf}{\nu_{\df}}
\newcommand{\nuF}{\nu_\tF}

\newcommand{\dNuF}{\nu_{\dtF}}
\newcommand{\ddNuF}{\nu_{\ddtF}}
\newcommand{\tilNu}{\tilde{\nu}}

\newcommand{\tilNuf}{\tilNu_f}
\newcommand{\tilNuF}{\tilNu_\tF}

\newcommand{\tilDuNuF}{\tilNu_{\dtF}}

\newcommand{\da}{\du{a}}

\newcommand{\dda}{\ddu{a}}

\newcommand{\dt}{{_\bullet}}

\newcommand{\brZ}{\bar{z}}

\newcommand{\nXz}{{\bf z}}
\newcommand{\nXx}{{\bf x}}
\newcommand{\nXy}{{\bf y}}

\newcommand{\tanBun}{\Lambda}
\newcommand{\coTanBun}{\tanBun^*}

\newcommand{\dInx}{\varpi}

\newcommand{\Val}{Val}

\newcommand{\iVal}{\Val^{-1}}

\newcommand{\Polygon}{\Delta}

\newcommand{\PolyHedA}{\sigma}
\newcommand{\PolyHedB}{\tau}
\newcommand{\curve}{\gamma}

\newcommand{\dEp}{\du{\ep}}
\newcommand{\ddEp}{\ddu{\ep}}

\newcommand{\lC}[1]{L_{#1}}
\newcommand{\epMat}{\trop{E}}
\newcommand{\dEpMat}{\du{\epMat}}
\newcommand{\ddEpMat}{\ddu{\epMat}}

\newcommand{\TrS}{\oplus}
\newcommand{\TrP}{\odot}

\newcommand{\maxPlus}{\Real \cup \{\tUniS\}}
\newcommand{\maxPlusSR}{(\maxPlus,\TrS,\TrP)}

\newcommand{\CH}{\mathcal{CH}}
\newcommand{\bCH}{\Game \CH}
\newcommand{\lbCH}{\underline{\bCH}}

\newcommand{\tropDet}[1]{TropDet(\trop{#1})}

\newcommand{\OP}{\left(}
\newcommand{\CP}{\right)}

\newcommand{\hlf}{\frac{1}{2}}
    \ifx\proof\undefined
    \newenvironment{proof}{
    \smallskip
    \noindent\emph{Proof.}}{\hfill\(\Box\)
    \bigskip
    } \fi

\newcommand{\printExtend}{tru}

\newcommand{\ifdef}[3]{\ifthenelse{\equal{#1}{true}}{#2}{#3}}

\newcommand {\secSpc} {\vskip 0cm}
\newcommand {\subSecSpc} {\vskip 0cm}

\newcommand {\parSpc} {\vskip 0.2cm}

\numberwithin{equation}{subsection}

%% file: diagram.tex
\makeatletter

\def\diagram{\m@th\leftwidth=\z@ \rightwidth=\z@ \topheight=\z@
\botheight=\z@ \setbox\@picbox\hbox\bgroup}

\def\enddiagram{\egroup\wd\@picbox\rightwidth\unitlength
\ht\@picbox\topheight\unitlength \dp\@picbox\botheight\unitlength
\hskip\leftwidth\unitlength\box\@picbox}

\def\bfig{\begin{diagram}}
\def\efig{\end{diagram}}
\newcount\wideness \newcount\leftwidth \newcount\rightwidth
\newcount\highness \newcount\topheight \newcount\botheight

\def\ratchet#1#2{\ifnum#1<#2 \global #1=#2 \fi}

\def\putbox(#1,#2)#3{%
\horsize{\wideness}{#3} \divide\wideness by 2
{\advance\wideness by #1 \ratchet{\rightwidth}{\wideness}}
{\advance\wideness by -#1 \ratchet{\leftwidth}{\wideness}}
\vertsize{\highness}{#3} \divide\highness by 2
{\advance\highness by #2 \ratchet{\topheight}{\highness}}
{\advance\highness by -#2 \ratchet{\botheight}{\highness}}
\put(#1,#2){\makebox(0,0){$#3$}}}

\def\putlbox(#1,#2)#3{%
\horsize{\wideness}{#3}
{\advance\wideness by #1 \ratchet{\rightwidth}{\wideness}}
{\ratchet{\leftwidth}{-#1}}
\vertsize{\highness}{#3} \divide\highness by 2
{\advance\highness by #2 \ratchet{\topheight}{\highness}}
{\advance\highness by -#2 \ratchet{\botheight}{\highness}}
\put(#1,#2){\makebox(0,0)[l]{$#3$}}}

\def\putrbox(#1,#2)#3{%
\horsize{\wideness}{#3}
{\ratchet{\rightwidth}{#1}}
{\advance\wideness by -#1 \ratchet{\leftwidth}{\wideness}}
\vertsize{\highness}{#3} \divide\highness by 2
{\advance\highness by #2 \ratchet{\topheight}{\highness}}
{\advance\highness by -#2 \ratchet{\botheight}{\highness}}
\put(#1,#2){\makebox(0,0)[r]{$#3$}}}

\def\adjust[#1]{} % For compatibility

\newcount \coefa
\newcount \coefb
\newcount \coefc
\newcount\tempcounta
\newcount\tempcountb
\newcount\tempcountc
\newcount\tempcountd
\newcount\xext
\newcount\yext
\newcount\xoff
\newcount\yoff
\newcount\gap%
\newcount\arrowtypea
\newcount\arrowtypeb
\newcount\arrowtypec
\newcount\arrowtyped
\newcount\arrowtypee
\newcount\height
\newcount\width
\newcount\xpos
\newcount\ypos
\newcount\run
\newcount\rise
\newcount\arrowlength
\newcount\halflength
\newcount\arrowtype
\newdimen\tempdimen
\newdimen\xlen
\newdimen\ylen
\newsavebox{\tempboxa}%
\newsavebox{\tempboxb}%
\newsavebox{\tempboxc}%

\newdimen\w@dth

\def\setw@dth#1#2{\setbox\z@\hbox{\m@th$#1$}\w@dth=\wd\z@
\setbox\@ne\hbox{\m@th$#2$}\ifnum\w@dth<\wd\@ne \w@dth=\wd\@ne \fi
\advance\w@dth by 1.2em}

\def\t@^#1_#2{\allowbreak\def\n@one{#1}\def\n@two{#2}\mathrel
{\setw@dth{#1}{#2}
\mathop{\hbox to \w@dth{\rightarrowfill}}\limits
\ifx\n@one\empty\else ^{\box\z@}\fi
\ifx\n@two\empty\else _{\box\@ne}\fi}}
\def\t@@^#1{\@ifnextchar_{\t@^{#1}}{\t@^{#1}_{}}}
\def\to{\@ifnextchar^{\t@@}{\t@@^{}}}

\def\t@left^#1_#2{\def\n@one{#1}\def\n@two{#2}\mathrel{\setw@dth{#1}{#2}
\mathop{\hbox to \w@dth{\leftarrowfill}}\limits
\ifx\n@one\empty\else ^{\box\z@}\fi
\ifx\n@two\empty\else _{\box\@ne}\fi}}
\def\t@@left^#1{\@ifnextchar_{\t@left^{#1}}{\t@left^{#1}_{}}}
\def\toleft{\@ifnextchar^{\t@@left}{\t@@left^{}}}

\def\two@^#1_#2{\allowbreak
\def\n@one{#1}\def\n@two{#2}\mathrel{\setw@dth{#1}{#2}
\mathop{\vcenter{\lineskip\z@\baselineskip\z@
                 \hbox to \w@dth{\rightarrowfill}%
                 \hbox to \w@dth{\rightarrowfill}}%
       }\limits
\ifx\n@one\empty\else ^{\box\z@}\fi
\ifx\n@two\empty\else _{\box\@ne}\fi}}
\def\tw@@^#1{\@ifnextchar _{\two@^{#1}}{\two@^{#1}_{}}}
\def\two{\@ifnextchar ^{\tw@@}{\tw@@^{}}}

\def\tofr@^#1_#2{\def\n@one{#1}\def\n@two{#2}\mathrel{\setw@dth{#1}{#2}
\mathop{\vcenter{\hbox to \w@dth{\rightarrowfill}\kern-1.7ex
                 \hbox to \w@dth{\leftarrowfill}}%
       }\limits
\ifx\n@one\empty\else ^{\box\z@}\fi
\ifx\n@two\empty\else _{\box\@ne}\fi}}
\def\t@fr@^#1{\@ifnextchar_ {\tofr@^{#1}}{\tofr@^{#1}_{}}}
\def\tofro{\@ifnextchar^ {\t@fr@}{\t@fr@^{}}}

\def\mon{\mathop{\m@th\hbox to
      14.6\P@{\lasyb\char'51\hskip-2.1\P@$\arrext$\hss
$\mathord\rightarrow$}}\limits} % width of \epi
\def\leftmono{\mathrel{\m@th\hbox to
14.6\P@{$\mathord\leftarrow$\hss$\arrext$\hskip-2.1\P@\lasyb\char'50%
}}\limits} % width of \epi
\mathchardef\arrext="0200       % amr minus for arrow extension (see \into)

\def\settypes(#1,#2,#3){\arrowtypea#1 \arrowtypeb#2 \arrowtypec#3}
\def\settoheight#1#2{\setbox\@tempboxa\hbox{#2}#1\ht\@tempboxa\relax}%
\def\settodepth#1#2{\setbox\@tempboxa\hbox{#2}#1\dp\@tempboxa\relax}%
\def\settokens`#1`#2`#3`#4`{%
     \def\tokena{#1}\def\tokenb{#2}\def\tokenc{#3}\def\tokend{#4}}
\def\setsqparms[#1`#2`#3`#4;#5`#6]{%
\arrowtypea #1
\arrowtypeb #2
\arrowtypec #3
\arrowtyped #4
\width #5
\height #6
}
\def\setpos(#1,#2){\xpos=#1 \ypos#2}

\def\settriparms[#1`#2`#3;#4]{\settripairparms[#1`#2`#3`1`1;#4]}%

\def\settripairparms[#1`#2`#3`#4`#5;#6]{%
\arrowtypea #1
\arrowtypeb #2
\arrowtypec #3
\arrowtyped #4
\arrowtypee #5
\width #6
\height #6
}

\def\resetparms{\settripairparms[1`1`1`1`1;500]\width 500}%default values%

\resetparms

\def\mvector(#1,#2)#3{%%
\put(0,0){\vector(#1,#2){#3}}%
\put(0,0){\vector(#1,#2){26}}%
}
\def\evector(#1,#2)#3{{%%
\arrowlength #3
\put(0,0){\vector(#1,#2){\arrowlength}}%
\advance \arrowlength by-30
\put(0,0){\vector(#1,#2){\arrowlength}}%
}}

\def\horsize#1#2{%
\settowidth{\tempdimen}{$#2$}%
#1=\tempdimen
\divide #1 by\unitlength
}

\def\vertsize#1#2{%
\settoheight{\tempdimen}{$#2$}%
#1=\tempdimen
\settodepth{\tempdimen}{$#2$}%
\advance #1 by\tempdimen
\divide #1 by\unitlength
}

\def\putvector(#1,#2)(#3,#4)#5#6{{%
\ifnum3<\arrowtype
\putdashvector(#1,#2)(#3,#4)#5\arrowtype
\else
\ifnum\arrowtype<-3
\putdashvector(#1,#2)(#3,#4)#5\arrowtype
\else
\xpos=#1
\ypos=#2
\run=#3
\rise=#4
\arrowlength=#5
\ifnum \arrowtype<0
    \ifnum \run=0
        \advance \ypos by-\arrowlength
    \else
        \tempcounta \arrowlength
        \multiply \tempcounta by\rise
        \divide \tempcounta by\run
        \ifnum\run>0
            \advance \xpos by\arrowlength
            \advance \ypos by\tempcounta
        \else
            \advance \xpos by-\arrowlength
            \advance \ypos by-\tempcounta
        \fi
    \fi
    \multiply \arrowtype by-1
    \multiply \rise by-1
    \multiply \run by-1
\fi
\ifcase \arrowtype
\or \put(\xpos,\ypos){\vector(\run,\rise){\arrowlength}}%
\or \put(\xpos,\ypos){\mvector(\run,\rise)\arrowlength}%
\or \put(\xpos,\ypos){\evector(\run,\rise){\arrowlength}}%
\fi\fi\fi
}}

\def\putsplitvector(#1,#2)#3#4{%%
\xpos #1
\ypos #2
\arrowtype #4
\halflength #3
\arrowlength #3
\gap 140
\advance \halflength by-\gap
\divide \halflength by2
\ifnum\arrowtype>0
   \ifcase \arrowtype
   \or \put(\xpos,\ypos){\line(0,-1){\halflength}}%
       \advance\ypos by-\halflength
       \advance\ypos by-\gap
       \put(\xpos,\ypos){\vector(0,-1){\halflength}}%
   \or \put(\xpos,\ypos){\line(0,-1)\halflength}%
       \put(\xpos,\ypos){\vector(0,-1)3}%
       \advance\ypos by-\halflength
       \advance\ypos by-\gap
       \put(\xpos,\ypos){\vector(0,-1){\halflength}}%
   \or \put(\xpos,\ypos){\line(0,-1)\halflength}%
       \advance\ypos by-\halflength
       \advance\ypos by-\gap
       \put(\xpos,\ypos){\evector(0,-1){\halflength}}%
   \fi
\else \arrowtype=-\arrowtype
   \ifcase\arrowtype
   \or \advance \ypos by-\arrowlength
       \put(\xpos,\ypos){\line(0,1){\halflength}}%
       \advance\ypos by\halflength
       \advance\ypos by\gap
       \put(\xpos,\ypos){\vector(0,1){\halflength}}%
   \or \advance \ypos by-\arrowlength
       \put(\xpos,\ypos){\line(0,1)\halflength}%
       \put(\xpos,\ypos){\vector(0,1)3}%
       \advance\ypos by\halflength
       \advance\ypos by\gap
       \put(\xpos,\ypos){\vector(0,1){\halflength}}%
   \or \advance \ypos by-\arrowlength
       \put(\xpos,\ypos){\line(0,1)\halflength}%
       \advance\ypos by\halflength
       \advance\ypos by\gap
       \put(\xpos,\ypos){\evector(0,1){\halflength}}%
   \fi
\fi
}

\def\putmorphism(#1)(#2,#3)[#4`#5`#6]#7#8#9{{%
\run #2
\rise #3
\ifnum\rise=0
  \puthmorphism(#1)[#4`#5`#6]{#7}{#8}#9%
\else\ifnum\run=0
  \putvmorphism(#1)[#4`#5`#6]{#7}{#8}#9%
\else
\setpos(#1)%
\arrowlength #7
\arrowtype #8
\ifnum\run=0
\else\ifnum\rise=0
\else
\ifnum\run>0
    \coefa=1
\else
   \coefa=-1
\fi
\ifnum\arrowtype>0
   \coefb=0
   \coefc=-1
\else
   \coefb=\coefa
   \coefc=1
   \arrowtype=-\arrowtype
\fi
\width=2
\multiply \width by\run
\divide \width by\rise
\ifnum \width<0  \width=-\width\fi
\advance\width by60
\if l#9 \width=-\width\fi
\putbox(\xpos,\ypos){#4}%            %node 1
{\multiply \coefa by\arrowlength%      %node 2
\advance\xpos by\coefa
\multiply \coefa by\rise
\divide \coefa by\run
\advance \ypos by\coefa
\putbox(\xpos,\ypos){#5} }%
{\multiply \coefa by\arrowlength%      %label
\divide \coefa by2
\advance \xpos by\coefa
\advance \xpos by\width
\multiply \coefa by\rise
\divide \coefa by\run
\advance \ypos by\coefa
\if l#9%
   \putrbox(\xpos,\ypos){#6}%
\else\if r#9%
   \putlbox(\xpos,\ypos){#6}%
\fi\fi }%
{\multiply \rise by-\coefc%             %arrow
\multiply \run by-\coefc
\multiply \coefb by\arrowlength
\advance \xpos by\coefb
\multiply \coefb by\rise
\divide \coefb by\run
\advance \ypos by\coefb
\multiply \coefc by70
\advance \ypos by\coefc
\multiply \coefc by\run
\divide \coefc by\rise
\advance \xpos by\coefc
\multiply \coefa by140
\multiply \coefa by\run
\divide \coefa by\rise
\advance \arrowlength by\coefa
\ifcase\arrowtype
\or \put(\xpos,\ypos){\vector(\run,\rise){\arrowlength}}%
\or \put(\xpos,\ypos){\mvector(\run,\rise){\arrowlength}}%
\or \put(\xpos,\ypos){\evector(\run,\rise){\arrowlength}}%
\fi}\fi\fi\fi\fi}}

\newcount\numbdashes \newcount\lengthdash \newcount\increment

\def\howmanydashes{% Actually returns both number and length
\numbdashes=\arrowlength \lengthdash=40
\divide\numbdashes by \lengthdash
\lengthdash=\arrowlength
\divide\lengthdash by \numbdashes
\increment=\lengthdash
\multiply\lengthdash by 3
\divide\lengthdash by 5
}

\def\putdashvector(#1)(#2,#3)#4#5{%
\ifnum#3=0 \putdashhvector(#1){#4}#5
\else
\ifnum#2=0
\putdashvvector(#1){#4}#5\fi\fi}

\def\putdashhvector(#1,#2)#3#4{{%
\arrowlength=#3 \howmanydashes
\multiput(#1,#2)(\increment,0){\numbdashes}%
{\vrule height .4pt width \lengthdash\unitlength}
\arrowtype=#4 \xpos=#1
\ifnum\arrowtype<0 \advance\arrowtype by 7 \fi
\ifcase\arrowtype
\or \advance\xpos by 10
    \put(\xpos,#2){\vector(-1,0){\lengthdash}}
    \advance\xpos by 40
    \put(\xpos,#2){\vector(-1,0){\lengthdash}}
\or \advance \xpos by 10
    \put(\xpos,#2){\vector(-1,0){\lengthdash}}
    \advance\xpos by  \arrowlength
    \advance\xpos by  -50
    \put(\xpos,#2){\vector(-1,0){\lengthdash}}
\or \advance\xpos by 10
    \put(\xpos,#2){\vector(-1,0){\lengthdash}}
\or \advance\xpos by \arrowlength
    \advance\xpos by -\lengthdash
    \put(\xpos,#2){\vector(1,0){\lengthdash}}
\or {\advance\xpos by 10
    \put(\xpos,#2){\vector(1,0){\lengthdash}}}
    \advance\xpos by \arrowlength
    \advance\xpos by -\lengthdash
    \put(\xpos,#2){\vector(1,0){\lengthdash}}
\or \advance\xpos by \arrowlength
    \advance\xpos by -\lengthdash
    \put(\xpos,#2){\vector(1,0){\lengthdash}}
    \advance\xpos by -40
    \put(\xpos,#2){\vector(1,0){\lengthdash}}
   \fi
}}

\def\putdashvvector(#1,#2)#3#4{{%
\arrowlength=#3 \howmanydashes
\ypos=#2 \advance\ypos by -\arrowlength
\multiput(#1,#2)(0,\increment){\numbdashes}%
    {\vrule width .4pt height \lengthdash\unitlength}
\arrowtype=#4 \ypos=#2
\ifnum\arrowtype<0 \advance\arrowtype by 7 \fi
\ifcase\arrowtype
\or \advance\ypos by \arrowlength \advance\ypos by -40
    \put(#1,\ypos){\vector(0,1){\lengthdash}}
    \advance\ypos by -40
    \put(#1,\ypos){\vector(0,1){\lengthdash}}
\or \advance\ypos by 10
    \put(#1,\ypos){\vector(0,1){\lengthdash}}
    \advance\ypos by \arrowlength \advance\ypos by -40
    \put(#1,\ypos){\vector(0,1){\lengthdash}}
\or \advance\ypos by \arrowlength \advance\ypos by -40
    \put(#1,\ypos){\vector(0,1){\lengthdash}}
\or \advance\ypos by 10
    \put(#1,\ypos){\vector(0,-1){\lengthdash}}
\or \advance\ypos by 10
    \put(#1,\ypos){\vector(0,-1){\lengthdash}}
    \advance\ypos by \arrowlength \advance\ypos by -40
    \put(#1,\ypos){\vector(0,-1){\lengthdash}}
\or \advance\ypos by 10
    \put(#1,\ypos){\vector(0,-1){\lengthdash}}
    \advance\ypos by 40
    \put(#1,\ypos){\vector(0,-1){\lengthdash}}
\fi
}}

\def\puthmorphism(#1,#2)[#3`#4`#5]#6#7#8{{%
\xpos #1
\ypos #2
\width #6
\arrowlength #6
\arrowtype=#7
\putbox(\xpos,\ypos){#3\vphantom{#4}}%
{\advance \xpos by\arrowlength
\putbox(\xpos,\ypos){\vphantom{#3}#4}}%
\horsize{\tempcounta}{#3}%
\horsize{\tempcountb}{#4}%
\divide \tempcounta by2
\divide \tempcountb by2
\advance \tempcounta by30
\advance \tempcountb by30
\advance \xpos by\tempcounta
\advance \arrowlength by-\tempcounta
\advance \arrowlength by-\tempcountb
\putvector(\xpos,\ypos)(1,0)\arrowlength\arrowtype
\divide \arrowlength by2
\advance \xpos by\arrowlength
\vertsize{\tempcounta}{#5}%
\divide\tempcounta by2
\advance \tempcounta by20
\if a#8 %
   \advance \ypos by\tempcounta
   \putbox(\xpos,\ypos){#5}%
\else
   \advance \ypos by-\tempcounta
   \putbox(\xpos,\ypos){#5}%
\fi}}

\def\putvmorphism(#1,#2)[#3`#4`#5]#6#7#8{{%
\xpos #1
\ypos #2
\arrowlength #6
\arrowtype #7
\settowidth{\xlen}{$#5$}%
\putbox(\xpos,\ypos){#3}%
{\advance \ypos by-\arrowlength
\putbox(\xpos,\ypos){#4}}%
{\advance\arrowlength by-140
\advance \ypos by-70
\ifdim\xlen>0pt
   \if m#8%
      \putsplitvector(\xpos,\ypos)\arrowlength\arrowtype
   \else
   \putvector(\xpos,\ypos)(0,-1)\arrowlength\arrowtype
   \fi
\else
   \putvector(\xpos,\ypos)(0,-1)\arrowlength\arrowtype
\fi}%
\ifdim\xlen>0pt
   \divide \arrowlength by2
   \advance\ypos by-\arrowlength
   \if l#8%
      \advance \xpos by-40
      \putrbox(\xpos,\ypos){#5}%
   \else\if r#8%
      \advance \xpos by40
      \putlbox(\xpos,\ypos){#5}%
   \else
      \putbox(\xpos,\ypos){#5}%
   \fi\fi
\fi
}}

\def\putsquarep<#1>(#2)[#3;#4`#5`#6`#7]{{%
\setsqparms[#1]%
\setpos(#2)%
\settokens`#3`%
\puthmorphism(\xpos,\ypos)[\tokenc`\tokend`{#7}]{\width}{\arrowtyped}b%
\advance\ypos by \height
\puthmorphism(\xpos,\ypos)[\tokena`\tokenb`{#4}]{\width}{\arrowtypea}a%
\putvmorphism(\xpos,\ypos)[``{#5}]{\height}{\arrowtypeb}l%
\advance\xpos by \width
\putvmorphism(\xpos,\ypos)[``{#6}]{\height}{\arrowtypec}r%
}}

\def\putsquare{\@ifnextchar <{\putsquarep}{\putsquarep%
   <\arrowtypea`\arrowtypeb`\arrowtypec`\arrowtyped;\width`\height>}}
\def\squaree{\@ifnextchar< {\squarep}{\squarep
   <\arrowtypea`\arrowtypeb`\arrowtypec`\arrowtyped;\width`\height>}}
\def\squarep<#1>[#2`#3`#4`#5;#6`#7`#8`#9]{{%       %     #2------>#3
\setsqparms[#1]%                                   %      |       |
\diagram%                                          %      |       |
\putsquarep<\arrowtypea`\arrowtypeb`\arrowtypec`%  %    #7|       |#8
\arrowtyped;\width`\height>%                       %      |       |
(0,0)[#2`#3`#4`{#5};#6`#7`#8`{#9}]%                %      |       |
\enddiagram%                                       %      v       v
}}                                                 %     #4------>#5
\def\putptrianglep<#1>(#2,#3)[#4`#5`#6;#7`#8`#9]{{%
\settriparms[#1]%
\xpos=#2 \ypos=#3
\advance\ypos by \height
\puthmorphism(\xpos,\ypos)[#4`#5`{#7}]{\height}{\arrowtypea}a%
\putvmorphism(\xpos,\ypos)[`#6`{#8}]{\height}{\arrowtypeb}l%
\advance\xpos by\height
\putmorphism(\xpos,\ypos)(-1,-1)[``{#9}]{\height}{\arrowtypec}r%
}}

\def\putptriangle{\@ifnextchar <{\putptrianglep}{\putptrianglep
   <\arrowtypea`\arrowtypeb`\arrowtypec;\height>}}
\def\ptriangle{\@ifnextchar <{\ptrianglep}{\ptrianglep
   <\arrowtypea`\arrowtypeb`\arrowtypec;\height>}}
\def\ptrianglep<#1>[#2`#3`#4;#5`#6`#7]{{%%    %      #2----->#3
\settriparms[#1]%                             %      |      /
\diagram%                                     %      |     /
\putptrianglep<\arrowtypea`\arrowtypeb`%      %    #6|    /#7
\arrowtypec;\height>%                         %      |   /
(0,0)[#2`#3`#4;#5`#6`{#7}]%                   %      |  /
\enddiagram%%                                 %      v v
}}                                            %      #4

\def\putqtrianglep<#1>(#2,#3)[#4`#5`#6;#7`#8`#9]{{%
\settriparms[#1]%
\xpos=#2 \ypos=#3
\advance\ypos by\height
\puthmorphism(\xpos,\ypos)[#4`#5`{#7}]{\height}{\arrowtypea}a%
\putmorphism(\xpos,\ypos)(1,-1)[``{#8}]{\height}{\arrowtypeb}l%
\advance\xpos by\height
\putvmorphism(\xpos,\ypos)[`#6`{#9}]{\height}{\arrowtypec}r%
}}

\def\putqtriangle{\@ifnextchar <{\putqtrianglep}{\putqtrianglep
   <\arrowtypea`\arrowtypeb`\arrowtypec;\height>}}
\def\qtriangle{\@ifnextchar <{\qtrianglep}{\qtrianglep
   <\arrowtypea`\arrowtypeb`\arrowtypec;\height>}}
\def\qtrianglep<#1>[#2`#3`#4;#5`#6`#7]{{%%    %        #2----->#3
\settriparms[#1]%                             %         \      |
\width=\height                                %          \     |
\diagram%                                     %         #6\    |#7
\putqtrianglep<\arrowtypea`\arrowtypeb`%      %            \   |
\arrowtypec;\height>%                         %             \  |
(0,0)[#2`#3`#4;#5`#6`{#7}]%                   %              v v
\enddiagram%%                                 %               #4
}}

\def\putdtrianglep<#1>(#2,#3)[#4`#5`#6;#7`#8`#9]{{%
\settriparms[#1]%
\xpos=#2 \ypos=#3
\puthmorphism(\xpos,\ypos)[#5`#6`{#9}]{\height}{\arrowtypec}b%
\advance\xpos by \height \advance\ypos by\height
\putmorphism(\xpos,\ypos)(-1,-1)[``{#7}]{\height}{\arrowtypea}l%
\putvmorphism(\xpos,\ypos)[#4``{#8}]{\height}{\arrowtypeb}r%
}}

\def\putdtriangle{\@ifnextchar <{\putdtrianglep}{\putdtrianglep
   <\arrowtypea`\arrowtypeb`\arrowtypec;\height>}}
\def\dtriangle{\@ifnextchar <{\dtrianglep}{\dtrianglep
   <\arrowtypea`\arrowtypeb`\arrowtypec;\height>}}
\def\dtrianglep<#1>[#2`#3`#4;#5`#6`#7]{{%%    %                  / |
\settriparms[#1]%                             %                 /  |
\width=\height                                %              #5/   |#6
\diagram%                                     %               /    |
\putdtrianglep<\arrowtypea`\arrowtypeb`%      %              /     |
\arrowtypec;\height>%                         %             v      v
(0,0)[#2`#3`#4;#5`#6`{#7}]%                   %            #3----->#4
\enddiagram%%                                 %                #7
}}

\def\putbtrianglep<#1>(#2,#3)[#4`#5`#6;#7`#8`#9]{{%
\settriparms[#1]%
\xpos=#2 \ypos=#3
\puthmorphism(\xpos,\ypos)[#5`#6`{#9}]{\height}{\arrowtypec}b%
\advance\ypos by\height
\putmorphism(\xpos,\ypos)(1,-1)[``{#8}]{\height}{\arrowtypeb}r%
\putvmorphism(\xpos,\ypos)[#4``{#7}]{\height}{\arrowtypea}l%
}}

\def\putbtriangle{\@ifnextchar <{\putbtrianglep}{\putbtrianglep
   <\arrowtypea`\arrowtypeb`\arrowtypec;\height>}}
\def\btriangle{\@ifnextchar <{\btrianglep}{\btrianglep
   <\arrowtypea`\arrowtypeb`\arrowtypec;\height>}}
\def\btrianglep<#1>[#2`#3`#4;#5`#6`#7]{{%%   %              | \
\settriparms[#1]%                            %              |  \
\width=\height                               %            #5|   \#6
\diagram%                                    %              |    \
\putbtrianglep<\arrowtypea`\arrowtypeb`%     %              |     \
\arrowtypec;\height>%                        %              v      v
(0,0)[#2`#3`#4;#5`#6`{#7}]%                  %              #3----->#4
\enddiagram%%                                %                 #7
}}

\def\putAtrianglep<#1>(#2,#3)[#4`#5`#6;#7`#8`#9]{{%
\settriparms[#1]%
\xpos=#2 \ypos=#3
{\multiply \height by2
\puthmorphism(\xpos,\ypos)[#5`#6`{#9}]{\height}{\arrowtypec}b}%
\advance\xpos by\height \advance\ypos by\height
\putmorphism(\xpos,\ypos)(-1,-1)[#4``{#7}]{\height}{\arrowtypea}l%
\putmorphism(\xpos,\ypos)(1,-1)[``{#8}]{\height}{\arrowtypeb}r%
}}

\def\putAtriangle{\@ifnextchar <{\putAtrianglep}{\putAtrianglep
   <\arrowtypea`\arrowtypeb`\arrowtypec;\height>}}
\def\Atriangle{\@ifnextchar <{\Atrianglep}{\Atrianglep
   <\arrowtypea`\arrowtypeb`\arrowtypec;\height>}}
\def\Atrianglep<#1>[#2`#3`#4;#5`#6`#7]{{%%         %         /   \
\settriparms[#1]%                                  %        /     \
\width=\height                                     %     #5/       \#6
\diagram%                                          %      /         \
\putAtrianglep<\arrowtypea`\arrowtypeb`%           %     /           \
\arrowtypec;\height>%                              %    v             v
(0,0)[#2`#3`#4;#5`#6`{#7}]%                        %   #3------------>#4
\enddiagram%%                                      %          #7
}}

\def\putAtrianglepairp<#1>(#2)[#3;#4`#5`#6`#7`#8]{{%
\settripairparms[#1]%
\setpos(#2)%
\settokens`#3`%
\puthmorphism(\xpos,\ypos)[\tokenb`\tokenc`{#7}]{\height}{\arrowtyped}b%
\advance\xpos by\height
\puthmorphism(\xpos,\ypos)[\phantom{\tokenc}`\tokend`{#8}]%
{\height}{\arrowtypee}b%
\advance\ypos by\height
\putmorphism(\xpos,\ypos)(-1,-1)[\tokena``{#4}]{\height}{\arrowtypea}l%
\putvmorphism(\xpos,\ypos)[``{#5}]{\height}{\arrowtypeb}m%
\putmorphism(\xpos,\ypos)(1,-1)[``{#6}]{\height}{\arrowtypec}r%
}}

\def\putAtrianglepair{\@ifnextchar <{\putAtrianglepairp}{\putAtrianglepairp%
   <\arrowtypea`\arrowtypeb`\arrowtypec`\arrowtyped`\arrowtypee;\height>}}
\def\Atrianglepair{\@ifnextchar <{\Atrianglepairp}{\Atrianglepairp%
   <\arrowtypea`\arrowtypeb`\arrowtypec`\arrowtyped`\arrowtypee;\height>}}

\def\Atrianglepairp<#1>[#2;#3`#4`#5`#6`#7]{{%           %  #2a
\settripairparms[#1]%                         %           / | \
\settokens`#2`%                               %          /  |  \
\width=\height                                %       #3/  #4   \#5
\diagram%                                     %        /    |    \
\putAtrianglepairp                            %       /     |     \
<\arrowtypea`\arrowtypeb`\arrowtypec`%        %      v      v      v
\arrowtyped`\arrowtypee;\height>%             %     #2b---->#2c---->#2d
(0,0)[{#2};#3`#4`#5`#6`{#7}]%                 %         #6     #7
\enddiagram%%
}}

\def\putVtrianglep<#1>(#2,#3)[#4`#5`#6;#7`#8`#9]{{%
\settriparms[#1]%
\xpos=#2 \ypos=#3
\advance\ypos by\height
{\multiply\height by2
\puthmorphism(\xpos,\ypos)[#4`#5`{#7}]{\height}{\arrowtypea}a}%
\putmorphism(\xpos,\ypos)(1,-1)[`#6`{#8}]{\height}{\arrowtypeb}l%
\advance\xpos by\height
\advance\xpos by\height
\putmorphism(\xpos,\ypos)(-1,-1)[``{#9}]{\height}{\arrowtypec}r%
}}

\def\putVtriangle{\@ifnextchar <{\putVtrianglep}{\putVtrianglep
   <\arrowtypea`\arrowtypeb`\arrowtypec;\height>}}
\def\Vtriangle{\@ifnextchar <{\Vtrianglep}{\Vtrianglep
   <\arrowtypea`\arrowtypeb`\arrowtypec;\height>}}
\def\Vtrianglep<#1>[#2`#3`#4;#5`#6`#7]{{%%     %        #2------------->#3
\settriparms[#1]%                              %         \             /
\width=\height                                 %          \           /
\diagram%                                      %         #6\         /#7
\putVtrianglep<\arrowtypea`\arrowtypeb`%       %            \       /
\arrowtypec;\height>%                          %             \     /
(0,0)[#2`#3`#4;#5`#6`{#7}]%                    %              v   v
\enddiagram%%                                  %               #4
}}

\def\putVtrianglepairp<#1>(#2)[#3;#4`#5`#6`#7`#8]{{
\settripairparms[#1]%
\setpos(#2)%
\settokens`#3`%
\advance\ypos by\height
\putmorphism(\xpos,\ypos)(1,-1)[`\tokend`{#6}]{\height}{\arrowtypec}l%
\puthmorphism(\xpos,\ypos)[\tokena`\tokenb`{#4}]{\height}{\arrowtypea}a%
\advance\xpos by\height
\puthmorphism(\xpos,\ypos)[\phantom{\tokenb}`\tokenc`{#5}]%
{\height}{\arrowtypeb}a%
\putvmorphism(\xpos,\ypos)[``{#7}]{\height}{\arrowtyped}m%
\advance\xpos by\height
\putmorphism(\xpos,\ypos)(-1,-1)[``{#8}]{\height}{\arrowtypee}r%
}}

\def\putVtrianglepair{\@ifnextchar <{\putVtrianglepairp}{\putVtrianglepairp%
    <\arrowtypea`\arrowtypeb`\arrowtypec`\arrowtyped`\arrowtypee;\height>}}
\def\Vtrianglepair{\@ifnextchar <{\Vtrianglepairp}{\Vtrianglepairp%
    <\arrowtypea`\arrowtypeb`\arrowtypec`\arrowtyped`\arrowtypee;\height>}}
\def\Vtrianglepairp<#1>[#2;#3`#4`#5`#6`#7]{{%  %  #2a---->#2b---->#2c
\settripairparms[#1]%                          %   \      |      /
\settokens`#2`%                                %    \     |     /
\diagram%                                      %   #5\   #6    /#7
\putVtrianglepairp                             %      \   |   /
<\arrowtypea`\arrowtypeb`\arrowtypec`%         %       \  |  /
\arrowtyped`\arrowtypee;\height>%              %        v v v
(0,0)[{#2};#3`#4`#5`#6`{#7}]%                  %         #2d
\enddiagram%%
}}

\def\putCtrianglep<#1>(#2,#3)[#4`#5`#6;#7`#8`#9]{{%
\settriparms[#1]%
\xpos=#2 \ypos=#3
\advance\ypos by\height
\putmorphism(\xpos,\ypos)(1,-1)[``{#9}]{\height}{\arrowtypec}l%
\advance\xpos by\height
\advance\ypos by\height
\putmorphism(\xpos,\ypos)(-1,-1)[#4`#5`{#7}]{\height}{\arrowtypea}l%
{\multiply\height by 2
\putvmorphism(\xpos,\ypos)[`#6`{#8}]{\height}{\arrowtypeb}r}%
}}

\def\putCtriangle{\@ifnextchar <{\putCtrianglep}{\putCtrianglep
    <\arrowtypea`\arrowtypeb`\arrowtypec;\height>}}
\def\Ctriangle{\@ifnextchar <{\Ctrianglep}{\Ctrianglep
    <\arrowtypea`\arrowtypeb`\arrowtypec;\height>}}
\def\Ctrianglep<#1>[#2`#3`#4;#5`#6`#7]{{%%   %                / |
\settriparms[#1]%                            %             #5/  |
\width=\height                               %              /   |
\diagram%                                    %             v    |
\putCtrianglep<\arrowtypea`\arrowtypeb`%     %           #3     |#6
\arrowtypec;\height>%                        %             \    |
(0,0)[#2`#3`#4;#5`#6`{#7}]%                  %            #7\   |
\enddiagram%%                                %               \  |
}}                                           %                v v
\def\putDtrianglep<#1>(#2,#3)[#4`#5`#6;#7`#8`#9]{{%
\settriparms[#1]%
\xpos=#2 \ypos=#3
\advance\xpos by\height \advance\ypos by\height
\putmorphism(\xpos,\ypos)(-1,-1)[``{#9}]{\height}{\arrowtypec}r%
\advance\xpos by-\height \advance\ypos by\height
\putmorphism(\xpos,\ypos)(1,-1)[`#5`{#8}]{\height}{\arrowtypeb}r%
{\multiply\height by 2
\putvmorphism(\xpos,\ypos)[#4`#6`{#7}]{\height}{\arrowtypea}l}%
}}

\def\putDtriangle{\@ifnextchar <{\putDtrianglep}{\putDtrianglep
    <\arrowtypea`\arrowtypeb`\arrowtypec;\height>}}
\def\Dtriangle{\@ifnextchar <{\Dtrianglep}{\Dtrianglep
   <\arrowtypea`\arrowtypeb`\arrowtypec;\height>}}
\def\Dtrianglep<#1>[#2`#3`#4;#5`#6`#7]{{%%  %          | \
\settriparms[#1]%                           %          |  \#6
\width=\height                              %          |   \
\diagram%                                   %          |    v
\putDtrianglep<\arrowtypea`\arrowtypeb`%    %        #5|    #3
\arrowtypec;\height>%                       %          |    /
(0,0)[#2`#3`#4;#5`#6`{#7}]%                 %          |   /#7
\enddiagram%%                               %          |  /
}}                                          %          v v
\def\setrecparms[#1`#2]{\width=#1 \height=#2}%
\def\recursep<#1`#2>[#3;#4`#5`#6`#7`#8]{{\m@th
\width=#1 \height=#2
\settokens`#3`
\settowidth{\tempdimen}{$\tokena$}
\ifdim\tempdimen=0pt
  \savebox{\tempboxa}{\hbox{$\tokenb$}}%
  \savebox{\tempboxb}{\hbox{$\tokend$}}%
  \savebox{\tempboxc}{\hbox{$#6$}}%
\else
  \savebox{\tempboxa}{\hbox{$\hbox{$\tokena$}\times\hbox{$\tokenb$}$}}%
  \savebox{\tempboxb}{\hbox{$\hbox{$\tokena$}\times\hbox{$\tokend$}$}}%
  \savebox{\tempboxc}{\hbox{$\hbox{$\tokena$}\times\hbox{$#6$}$}}%
\fi
\ypos=\height
\divide\ypos by 2
\xpos=\ypos
\advance\xpos by \width
\bfig
\putCtrianglep<-1`1`1;\ypos>(0,0)[`\tokenc`;#5`#6`{#7}]%
\puthmorphism(\ypos,0)[\tokend`\usebox{\tempboxb}`{#8}]{\width}{-1}b%
\puthmorphism(\ypos,\height)[\tokenb`\usebox{\tempboxa}`{#4}]{\width}{-1}a%
\advance\ypos by \width
\putvmorphism(\ypos,\height)[``\usebox{\tempboxc}]{\height}1r%
\efig
}}

\def\recurse{\@ifnextchar <{\recursep}{\recursep<\width`\height>}}

\def\puttwohmorphisms(#1,#2)[#3`#4;#5`#6]#7#8#9{{%
\puthmorphism(#1,#2)[#3`#4`]{#7}0a
\ypos=#2
\advance\ypos by 20
\puthmorphism(#1,\ypos)[\phantom{#3}`\phantom{#4}`#5]{#7}{#8}a
\advance\ypos by -40
\puthmorphism(#1,\ypos)[\phantom{#3}`\phantom{#4}`#6]{#7}{#9}b
}}

\def\puttwovmorphisms(#1,#2)[#3`#4;#5`#6]#7#8#9{{%
\putvmorphism(#1,#2)[#3`#4`]{#7}0a
\xpos=#1
\advance\xpos by -20
\putvmorphism(\xpos,#2)[\phantom{#3}`\phantom{#4}`#5]{#7}{#8}l
\advance\xpos by 40
\putvmorphism(\xpos,#2)[\phantom{#3}`\phantom{#4}`#6]{#7}{#9}r
}}

\def\puthcoequalizer(#1)[#2`#3`#4;#5`#6`#7]#8#9{{%
\setpos(#1)%
\puttwohmorphisms(\xpos,\ypos)[#2`#3;#5`#6]{#8}11%
\advance\xpos by #8
\puthmorphism(\xpos,\ypos)[\phantom{#3}`#4`#7]{#8}1{#9}
}}

\def\putvcoequalizer(#1)[#2`#3`#4;#5`#6`#7]#8#9{{%
\setpos(#1)%
\puttwovmorphisms(\xpos,\ypos)[#2`#3;#5`#6]{#8}11%
\advance\ypos by -#8
\putvmorphism(\xpos,\ypos)[\phantom{#3}`#4`#7]{#8}1{#9}
}}

\def\putthreehmorphisms(#1)[#2`#3;#4`#5`#6]#7(#8)#9{{%
\setpos(#1) \settypes(#8)
\if a#9 %
     \vertsize{\tempcounta}{#5}%
     \vertsize{\tempcountb}{#6}%
     \ifnum \tempcounta<\tempcountb \tempcounta=\tempcountb \fi
\else
     \vertsize{\tempcounta}{#4}%
     \vertsize{\tempcountb}{#5}%
     \ifnum \tempcounta<\tempcountb \tempcounta=\tempcountb \fi
\fi
\advance \tempcounta by 60
\puthmorphism(\xpos,\ypos)[#2`#3`#5]{#7}{\arrowtypeb}{#9}
\advance\ypos by \tempcounta
\puthmorphism(\xpos,\ypos)[\phantom{#2}`\phantom{#3}`#4]{#7}{\arrowtypea}{#9}
\advance\ypos by -\tempcounta \advance\ypos by -\tempcounta
\puthmorphism(\xpos,\ypos)[\phantom{#2}`\phantom{#3}`#6]{#7}{\arrowtypec}{#9}
}}

\def\setarrowtoks[#1`#2`#3`#4`#5`#6]{%
\def\toka{#1}
\def\tokb{#2}
\def\tokc{#3}
\def\tokd{#4}
\def\toke{#5}
\def\tokf{#6}
}
\def\hex{\@ifnextchar <{\hexp}{\hexp<1000`400>}}
\def\hexp<#1`#2>[#3`#4`#5`#6`#7`#8;#9]{%
\setarrowtoks[#9]
\yext=#2 \advance \yext by #2
\xext=#1 \advance\xext by \yext
\bfig
\putCtriangle<-1`0`1;#2>(0,0)[`#5`;\tokb``\tokd]
\xext=#1 \yext=#2 \advance \yext by #2
\putsquare<1`0`0`1;\xext`\yext>(#2,0)[#3`#4`#7`#8;\toka```\tokf]
\advance \xext by #2
\putDtriangle<0`1`-1;#2>(\xext,0)[`#6`;`\tokc`\toke]
\efig
}
\makeatother